\documentclass[a4paper,11pt,reqno]{amsart}
\usepackage{color, graphicx, enumerate, amssymb,bbm}
\usepackage{hyperref}
\usepackage{epsfig,wrapfig}
\usepackage{pxfonts}
\usepackage{graphicx}

\usepackage{eucal}

\usepackage[utf8]{inputenc}

\usepackage{amsmath}

\usepackage{geometry}
\usepackage{enumitem}

\newtheorem{theorem}{Theorem}[section] 
\newtheorem{lemma}[theorem]{Lemma}

\newtheorem{definition}{Definition}[section]
\newtheorem{proposition}[theorem]{Proposition}
\newtheorem{remark}[theorem]{Remark}

\newcommand{\bR}{{\mathbb R}}

\newcommand{\inc}[1]{\hyperref[inc]{{\normalfont(inc){\ensuremath{_{#1}}}}}}
\newcommand{\dec}[1]{\hyperref[dec]{{\normalfont(dec){\ensuremath{_{#1}}}}}}
\newcommand{\azero}{\hyperref[azero]{{\normalfont(A0)}}}
\newcommand{\auno}{\hyperref[azero]{{\normalfont(A1)}}}
\newcommand{\vauno}{\hyperref[va1]{{\normalfont(VA1)}}}

\newcommand{\Phiw}{\Phi_{\rm{w}}}
\newcommand{\Phic}{\Phi_{\rm{c}}}

 \providecommand{\Xint}[1]{\mathchoice
    {\XXint\displaystyle\textstyle{#1}}%
    {\XXint\textstyle\scriptstyle{#1}}%
    {\XXint\scriptstyle\scriptscriptstyle{#1}}%
    {\XXint\scriptscriptstyle\scriptscriptstyle{#1}}%
    \!\int}
  \providecommand{\XXint}[3]{{\setbox0=\hbox{$#1{#2#3}{\int}$}
      \vcenter{\hbox{$#2#3$}}\kern-.5\wd0}}
  
  \providecommand{\dashint}{\mathop{\Xint-}}
  \definecolor{MyDarkGreen}{rgb}{0,0.50,0.04}

\numberwithin{equation}{section}

\title[$\varphi$ almost-minimizers]{Lipschitz regularity of almost-minimizers in one-phase problems with generalized Orlicz growth}

\author{Chiara Leone}
\address[Chiara Leone]{Department of Mathematics and Applications ``R. Caccioppoli'', University of Naples Federico II, Via Cintia, Monte S. Angelo, 80126 Naples, Italy}
\email{chiara.leone@unina.it}
 
\author{Giovanni Scilla}
\address[Giovanni Scilla]{Department of Mathematics and Applications ``R. Caccioppoli'', University of Naples Federico II, Via Cintia, Monte S. Angelo, 80126 Naples, Italy}
\email[Giovanni Scilla]{giovanni.scilla@unina.it}

\author{Francesco Solombrino}
\address[Francesco Solombrino]{Department of Biological and Environmental Sciences and Technologies, University of Salento, Complesso Ecotekne, via Monteroni 165, 73100 Lecce, Italy}
\email{francesco.solombrino@unisalento.it}

\author{Anna Verde}
\address[Anna Verde]{Department of Mathematics and Applications ``R. Caccioppoli'', University of Naples Federico II, Via Cintia, Monte S. Angelo, 80126 Naples, Italy}
\email{anna.verde@unina.it}

\date{\today}

\begin{document}

\begin{abstract}

Optimal local Lipschitz regularity %For a given constant $\lambda > 0$ and a bounded Lipschitz domain $\Omega \subset \mathbb{R}^d$ ($d \geq 2$), we establish that 
for scalar almost minimizers of Alt-Caffarelli-type functionals
$$
 \mathcal{F}({v}; \Omega) = \int_\Omega \varphi(x,\left|\nabla v(x) \right|)+ \lambda \chi_{\{{v} >0\}} (x) \, \mathrm{d}x\,,
 $$
with growth function $\varphi$ a generalized Orlicz function, is established. %The boundary regularity result is new also in the case of $p$-growth, where the subquadratic regime $1<p<2$ was missing.}
\end{abstract}



\maketitle

\tableofcontents

\section{Introduction and the main results}
Our paper is concerned with giving a unifying perspective on the regularity theory of scalar-valued almost minimizers of an Alt--Caffarelli type functional with generalized Orlicz growth in the energy:
\begin{equation}\label{eq: model}
\int_\Omega \varphi(x, |\nabla u|) + \lambda \chi_{\{u>0\}} \, \mathrm{d}x\,,
\end{equation}
where \(\varphi\) models non-standard growth conditions, encompassing classical \(p\)-growth and Orlicz growth settings and allowing for inhomegeneities in the space variable.

The Alt--Caffarelli functional arises naturally in a variety of contexts involving phase separation phenomena, such as fluid interfaces, combustion, and optimal material design. Its central feature is the appearance of a free boundary -- i.e., the boundary of the region where the solution is positive -- which is not known a priori and must be determined as part of the problem. From a mathematical perspective, this functional serves as a canonical model in the study of free boundary problems. It brings together variational methods, geometric measure theory, and PDE analysis. Understanding the regularity of minimizers and their associated free boundaries is crucial not only for theoretical reasons but also for the stability and accuracy of computational methods in applications.

Nonlinear growth conditions induce additional complexities: most nonlinear regularity approaches must distinguish between subquadratic ($1<p<2$) and superquadratic ($p\geq2$) growth regimes. Structural assumptions on \(\varphi\), such as uniform convexity and appropriate control of the growth exponents, are essential to obtain Lipschitz regularity results. %Moreover, regularity up to the boundary is currently well-understood only in the superquadratic regime, as shown in \cite{BFS24}.

The Lipschitz continuity of minimizers to the Alt--Caffarelli functional is a fundamental property that guarantees the well-posedness of the problem and serves as a stepping stone toward free boundary regularity.

In the variational framework, minimizers (and more generally, almost minimizers) arise as solutions to energy minimization problems. The pioneering work of Alt and Caffarelli \cite{AltCaf} established Lipschitz regularity for scalar minimizers, leveraging a delicate blow-up analysis and monotonicity formulas. Caffarelli later refined these techniques \cite{Caff87, Caff88, Caff89}, introducing tools like Harnack-type inequalities and improvement of flatness. In its most general version, developed in an Orlicz-Sobolev setting in \cite{MW08}, 
this approach reformulates the minimization problem for $\mathcal{J}_{\mathrm{\varphi}, \lambda}(u,\Omega) := \int_\Omega \left(\mathrm{\varphi}(|\nabla u|)+\lambda\chi_{\{u>0\}}\right)\,\mathrm{d}x$ as a nonlinear PDE
$$
%\mathcal{J}_{\mathrm{G}, \lambda}(u,\Omega) := \int_\Omega \left(\mathrm{G}(|\nabla u|)+\lambda\chi_{\{u>0\}}\right)\,\mathrm{d}x  \quad \Rightarrow \quad 
\left\{\begin{array}{ccl}
    \mathrm{div} \left(\varphi'(|\nabla u|)\dfrac{\nabla u}{|\nabla u|}\right) = 0 & \text{in}&  \{u>0\}\cap \Omega \\
    u=0,\,\,\,\,\, |\nabla u| = \lambda^{\ast} & \text{on} & \partial \{u>0\}\cap \Omega\\
    \varphi'(\lambda^{\ast})\lambda^{\ast}-\varphi(\lambda^{\ast}) = \lambda%G^{\prime}(t) = g(t)& \text{and} & 
\end{array}
\right.
$$
with suitable boundary conditions. This weak formulation, and a sub- and supersolution method (in a suitable viscosity sense)  allow the authors to prove that solutions to the optimization problem are locally Lipschitz continuous. This extends the Alt-Caffarelli's results for the scenario of the Orlicz-Sobolev framework. Additionally, they address Caffarelli’s classification scheme: flat and Lipschitz free boundaries are locally $C^{1,\alpha}$ for some $\alpha(\verb"universal") \in (0,  1)$. Although not directly addressed in our contribution, let us also remark that a flourishing literature on the regularity of the free boundary has developed from these results,
particularly in the scalar case. Key contributions include \cite{AltCaf, Caff87, Caff88, Caff89, DephiSpoVel21, Velichkov23} for the scalar case, and \cite{ASUW15, Fotouhi2021, DeSJeonShah22, DeSJeonShah23, DeSJeonShah25} for the vectorial setting. 
%\begin{itemize}
%  \item Scalar case: 
%  \item Vectorial case: \cite{ASUW15, Fotouhi2021, DeSJeonShah22, DeSJeonShah23, DeSJeonShah25}
%\end{itemize}

The above tools are however not available when dealing with  almost minimizers. This concept, introduced in works such as De Silva and Savin \cite{DeSS20}, allows for perturbative and non-exact minimization, providing a robust framework to study solutions in inhomogeneous or approximate settings. Related contributions by David, Engelstein, Smit Vega and Toro \cite{DEST21,DET19, DT15} use compactness and approximation arguments to extend regularity results beyond the exact minimizer class. 

A remarkable approach to regularity of almost minimizers in a nonlinear $p$-Laplace setting has been later devised in \cite{DiPFFV24}. Their approach is based on local regularity estimates on the $p$-harmonic replacement of an almost minimizer and a dichotomy theory according to which, roughly speaking, the average of the energy of an almost minimizer decreases in a smaller ball, unless one is arbitrarily close to the case of linear functions. This analysis, concerning interior regularity, introduces some restrictions on the exponent $p$, and is essentially carried out with different techniques and estimates for the sub- and superquadratic case.

In the vectorial setting, Lipschitz regularity becomes even more subtle due to the interaction between components. Recent results by Bayrami, Fotouhi, and Shahgholian \cite{BFS24}, and by De Silva, Jeon, and Shahgholian \cite{DeSJeonShah22, DeSJeonShah23}, address weakly coupled systems and singular behaviors, obtaining Lipschitz bounds under structural and coupling assumptions. In particular, the paper \cite{BFS24} handles regularity (up to the boundary), obtained via a blow-up approach which will be also expedient for our analysis, as described below. 
%However, while interior regularity is obtained for an arbitrary exponent in the $p$-Laplace energy, the boundary regularity result is limited to the case $p\geq 2$.

The analysis has been further extended to the case of non-standard growth (e.g., Orlicz spaces) further enriching the theory of Lipschitz regularity of almost minimizers in non-classical environments. This has been mostly done for autonomous integrands (i.e., the case where $\varphi(x, t)=\varphi(t)$). The scalar case has been addressed  in \cite{dSSV}, and later extended to the vector valued case  in the recent contribution \cite{PdSY},  always dealing with a weakly coupled  system. We may also refer the reader to the introduction in \cite{PdSY} for a thorough analysis of the state-of-the-art for this kind of problems, as well as for a rich bibliography. \\

\noindent \emph{Description of our result.} We now come to the description of our results. In this paper, we prove Lipschitz regularity  for scalar-valued almost minimizers of \eqref{eq: model} by exploiting the very general  theory for non-autonomous functionals with Uhlenbeck structure introduced in \cite{HastoOk}. It is based on a single condition involving both the $x$ and $t$ directions (see Definition \ref{eleuteri} below), which allows one for comparison estimates with minimizers of a locally well-defined autonomous integrand $\widetilde \varphi$ whose properties are recalled in Section \ref{sec:regularizedphi}.

This approach can be succesfully combined with the scheme pursued in \cite{BFS24}. 
%entailing both interior and then boundary regularity results. 
The path we follow runs 
%indeed on parallel tracks, provided  one is able to replace interior with boundary regularity results for elliptic functionals, 
along these lines:

\begin{itemize}
\item at first, higher integrability estimates, and a reverse H\" older-type inequality for almost minimizers are established;
\item a local comparison estimate with the $\widetilde \varphi$-harmonic replacement allows one to recover a  Morrey-type estimate (see \eqref{eq:intMorrey}) and eventually $C^{\alpha}$ continuity of almost minimizers for any $\alpha$ in $(0,1)$;
\item a $C^{1, \alpha}$  regularity result {\it away from the free boundary} for suitable $\alpha$ is established by a compactness and lower semicontinuity argument on blown-up sequences of almost minimizers, which are shown to tend to a minimizer of an autonomous functional. Here, assumption \vauno{} in Definition \ref{eleuteri} plays a crucial role;
\item a key step, as in \cite{BFS24}, is finally to show that a bounded almost minimizer of \eqref{eq: model} is sublinear in a neighborhood of a free-boundary point. For this, in particular, the $C^{\alpha}$ continuity of almost minimizers is  exploited.

\end{itemize}

Observe that the combination of the estimates in \cite{HastoOk} with blow-up arguments causes some additional nontrivial difficulties to our analysis, as we must ensure that some constants, which depend themselves on the chosen local almost minimizer $u$, can also be used for providing uniform estimates for the blown-up sequences: this is apparent, for instance, in the proofs of  Propositions \ref{lem:scaledv} and \ref{prop:stimaM}.

The unified treatment we develop encompasses all the relevant examples in literature, provided the bulk energy has a growth from below with exponent $p>1$ and needs not distinguish between sub- and superquadratic energies. As relevant examples of energies undergoing non-standard growth we report here the perturbed {\em Orlicz}, the so-called {\em variable exponent}, 
%(~\cite{czech})
and the {\em double-phase} case
	\begin{equation*}
	\label{e:intro1}
	a(x)\varphi(|\xi|),\qquad
	|\xi|^{p(x)}, \qquad\text{and} \qquad |\xi|^{p} + a(x) |\xi|^{q} \qquad \text{for $(x, \xi) \in \mathbb{R}^{d} \times \mathbb{R}^{d}$,}
	\end{equation*}
%for suitable choices of the $N$-function $\varphi:[0,+\infty)\to[0,+\infty)$, of the exponent function $p \colon \R^{d} \to (1, +\infty)$, of the exponents  $1 < p < q <+\infty$, and of the weight function~$a \colon \R^{d} \to [0, +\infty)$, 
while an exhaustive list of  examples  can be found in \cite{HastoOk}.  
 %In particular, Lipschitz continuity up to the boundary is, to best of our knowledge, a novel result even in the case
%$\varphi(x, t)= t^p$ with $1<p<2$.

We remark that in this paper we limit ourselves to the case of scalar-valued almost minimizers but the same result can be
extended to the vectorial problem of a weakly coupled system with some additional, but manageable, effort, using the same procedure of \cite{BFS24}.
%essentially in order to exploit the methods of \cite{HastoOk} which are so far available only in that setting. 
%While regularity for a fairly general class of autonomous functionals with Orlicz growth has been very recently established in the vectorial setting in \cite{PdSY},  
 Extensions to both the autonomous and non-autonomous  vectorial case for strongly coupled systems represents a challenging direction of research which we plan to address in future contributions.

In order to introduce the main result of our paper, we formulate our problem and specify the definition of almost-minimizer.

Given a constant $\lambda>0$, and a bounded open set $\Omega \subset \bR^d$ ($d \geq 2$), we will deal with local almost-minimizers of the  functional
\begin{equation}
\label{E0}
  \mathcal{F}({v}; \Omega) := \int_\Omega \varphi(x,\left|\nabla v(x) \right|)+ \lambda \chi_{\{{v} >0\}} (x) \, \mathrm{d}x\,,
\end{equation}
for $v\in W^{1,\varphi}(\Omega)$ with $v\ge 0$.
%over  an admissible class 
%$$ \mathcal{A}:=\left\{{v} \in W^{1,\varphi}(\Omega) \, : \, \text{${v}={w}$ on $\partial \Omega$ and $v \geq 0$} \right\}, $$
%for a non-negative function $w\in W^{1,\varphi}(\Omega)$. 
The precise notion of almost minimizers that we use is the following.

\begin{definition}
We say that $u:\Omega \to \mathbb{R}$ is a (local) \emph{almost-minimizer} for $ \mathcal{F}$ in $\Omega$, with constant $\kappa$ and exponent $\beta$, if 
$$
 \mathcal{F}(u; B_r(x_0)) \leq \left(1+\kappa r^{\beta} \right)  \mathcal{F}(w; B_r(x_0)), 
$$
for every ball $B_r(x_0)$ such that $\overline{B_r(x_0)} \subset \Omega$ and every ${w} \in W^{1,\varphi}(B_r(x_0))$ such that ${u}={w}$ on $\partial B_r(x_0)$.
\end{definition}

The main result of the paper is the following. 

\begin{theorem}[Interior regularity]
\label{T1}
Let $\Omega\subset\mathbb{R}^d$ be a bounded open set. Let $\varphi \in \Phic(\Omega)$, $\varphi(x,\cdot)\in C^1([0,\infty))$ be satisfying \vauno{}, and such that $\varphi_t$ comply with \azero{}, \inc{p-1} and \dec{q-1} for some $1<p\leq q$. Let ${u} :\Omega \to \mathbb{R} $ be an almost-minimizer of $\mathcal{F}$ in $\Omega$. Then, ${u}$ is locally Lipschitz continuous in $\Omega$.
\end{theorem}

We are investigating the boundary version of the previous result in a paper currently in preparation.

\medskip

\emph{Outline of the paper.} The rest of the paper is organized as follows. In Section~\ref{sec: prel} we fix the basic notation and recall some basic facts and technical results about Orlicz and generalized Orlicz functions, together with some technical lemmas. Section \ref{sec:regauto} contains some supporting regularity results for autonomous problems in divergence form, which are exploited in Section~\ref{sec:compar}, where we obtain Caccioppoli type estimates and higher integrability results for almost minimizers, together with useful comparison estimates with the solution of a suitable autonomous problem. Section~\ref{sec:locholder} collects two main ingredients in order to get the main result: the local H\"{o}lder continuity of almost minimizers, Theorem~\ref{thm:C0alpha}, and that of their gradients away from the free boundary, Theorem~ \ref{thm:c1alphabound}. Section~\ref{sec:lipschitzcontinuous} is entirely devoted to the proof of the main result: the main steps are Lemma~\ref{lem:boundgrad0}, where an interior uniform bound for the gradient of an almost minimizer is provided; Proposition~\ref{lem:scaledv}, showing that a suitable blow-up sequence of almost minimizers converges to the solution of a limit  autonomous problem, and its consequence Proposition~\ref{prop:stimaM}, where the sublinearity of a bounded almost minimizer in a neighborhood of a free-boundary point is shown. Finally, in Appendix~\ref{sec:appendix}, we collects some technical results mainly employed in the proof of Lemma \ref{lem:asymptphi} and Proposition~\ref{lem:scaledv}.

\section{Basic notation and preliminaries}\label{sec: prel}

 We start with some basic notation.   Let $\Omega \subset \mathbb{R}^d$  be  open and bounded.
For every $x\in \mathbb{R}^d$ and $r>0$ we indicate by $B_r(x) \subset \mathbb{R}^d$ the open ball with center $x$ and radius $r$. We will often use the shorthand \(B_r\) when either $x=0$ or the center $x$ is not relevant.  For $x$, $y\in \mathbb{R}^d$, we use the notation $x\cdot y$ for the scalar product and $|x|$ for the  Euclidean  norm.  The $m$-dimensional Lebesgue measure of the unit ball in $\mathbb{R}^m$ is indicated by $\gamma_m$ for every $m \in \mathbb{N}$.   We denote by $\mathcal{L}^d$  the $d$-dimensional Lebesgue measure.  %For $A \subset \R^d$, $
The closure of $A$ is denoted by $\overline{A}$. The diameter of $A$ is indicated by ${\rm diam}(A)$. We write $\chi_A$ for the  characteristic  function of any $A\subset  \mathbb{R}^d$, which is 1 on $A$ and 0 otherwise.  

Given two functions $f,g:[0,+\infty)\to\mathbb{R}$, we write $f\sim g$, and we say that $f$ and $g$ are equivalent, if there exist constants $c_{1}, c_{2} >0$ such that $c_{1}g(t) \leq f(t) \leq c_{2}g(t)$ for any $t\ge 0$. Similarly the symbol $\lesssim$ stands for $\le$ up to a constant. {$L^0(\Omega)$ denotes the set of the measurable functions on $\Omega$.}

\subsection{Generalized $\Phi$-functions and Orlicz spaces}\label{sec:genorl}

We introduce some basic definitions and useful facts about generalized $\Phi$-functions and Orlicz spaces, only considering concepts we will use. We refer the reader to \cite{HH} for a comprehensive treatment of the topic. 

{\begin{definition}\label{weakphi}
Let $\varphi: [0,+\infty)\to [0, +\infty]$ be increasing with $\varphi(0) = 0$, $\displaystyle\lim_{t\to 0^+} \varphi(t) = 0$ and $\displaystyle\lim_{t\to+\infty}\varphi(t) = +\infty$. Such $\varphi$ is called a
\begin{enumerate}
\item[(i)] \emph{weak $\Phi$-function} if $\frac{\varphi(t)}{t}$ is \emph{almost increasing}, meaning that there exists $L\ge 1$ such that $\frac{\varphi(t)}{t}\le L\frac{\varphi(s)}{s}$ for $0<t\le s$.
\item[(ii)] \emph{convex $\Phi$-function} if $\varphi$ is left-continuous and convex. 
\end{enumerate}
\end{definition}
By virtue of Remark \ref{rem:varinc}, each convex $\Phi$-function is a weak $\Phi$-function. If $\varphi$ is a convex $\Phi$-function, then there exists $\varphi'$ the \emph{right derivative} of $\varphi$, which is non-decreasing and right-continuous, and such that
\begin{equation*}
\varphi(t)=\int_0^t \varphi'(s)\,\mathrm{d}s \,.
\end{equation*}
A special subclass of convex $\Phi$-functions is represented by the %so called ``nice Young functions'', also known as 
$N$-functions (see, e.g., \cite[Ch.I]{KR}).}

{\begin{definition}\label{Nfunc}
A function $\varphi: [0, \infty)\rightarrow [0, \infty)$ is said to be an \emph{$N$-function} if it admits the representation
\begin{equation*}
\varphi(t)=\int_0^t a(\tau)\,\mathrm{d}\tau
\end{equation*}
where $a(s)$ is right-continuous, non-decreasing for $s>0$, $a(s)>0$ for $s>0$ and satisfies the conditions
\begin{equation}
a(0)=0\,,\quad \lim_{s\to+\infty} a(s) =+\infty\,.
\label{eq:asympright}
\end{equation}
%$\varphi(0)=0$ and 
%there exists a right continuous nondecreasing derivative $\varphi'$ satisfying $\varphi'(0)=0$,
%$\varphi'(t)>0$ for $t>0$ and $\displaystyle{\lim_{t\rightarrow \infty} \varphi'(t)=\infty}$. Especially $\varphi$ is convex. 
\end{definition}
The function $a(t)$ is nothing else than the right-derivative of $\varphi(t)$. %, which will be denoted by $\varphi'(t)$. 
As a straightforward consequence of the definition, we have that an $N$-function $\varphi$ is continuous, $\varphi(0)=0$ and $\varphi$ is increasing. Moreover, $\varphi$ is a convex function, and, in view of Remark \ref{rem:varinc}, it satisfies \inc{1}.}
{ Conditions \eqref{eq:asympright} imply
\begin{equation}
\lim_{t\to0^+} \frac{\varphi(t)}{t}=0\,,\quad \lim_{t\to +\infty} \frac{\varphi(t)}{t}=+\infty\,.
\label{eq:asympright2}
\end{equation}
It can be shown that an equivalent definition of $N$-function is the following: a continuous convex function $\varphi$ is called an $N$-function if it satisfies \eqref{eq:asympright2}. }

{For our purposes, we need functions $\varphi$ to depend also on the spatial variable $x$. }
{\begin{definition}\label{generalizedfunct}
Let $\varphi:\Omega\times[0,\infty)\to [0,\infty]$. We call 
$\varphi$ a \textit{generalized} weak $\Phi$-function (resp., convex $\Phi$-function, $N$-function) if
\begin{enumerate}
\item[(1)] $x\mapsto \varphi(x,|f(x)|)$ is measurable for every $f\in L^0(\Omega)$; 
\item[(2)] $t\mapsto \varphi(x,t)$ is a weak $\Phi$-function (resp., a convex $\Phi$-function, an $N$-function) for every $x\in\Omega$. 
\end{enumerate}
We write $\varphi\in\Phiw(\Omega)$, $\varphi\in\Phic(\Omega)$ and $\varphi\in N(\Omega)$, respectively. If $\varphi$ does not depend on $x$, we will adopt the shorthands $\varphi\in\Phiw$, $\varphi\in\Phic$ and $\varphi\in N$, respectively. For the right-derivative of a generalized convex $\Phi$-function, we will use the notation $\varphi_t$ in place of $\varphi'$.
\end{definition} }

{For a bounded function $\varphi:\Omega\times[0,+\infty)\to[0,+\infty)$ and a ball $B_r(x_0)\subset\Omega$ we define, for every $t\ge 0$,
\begin{equation}\label{phimeno}
\varphi_{r,x_0}^-(t):=\inf_{x\in B_r(x_0)}\varphi(x,t)\quad\hbox{ and }\quad \varphi_{r,x_0}^+(t) :=  \sup_{x\in B_r(x_0)}\varphi(x,t).
\end{equation}}

%Let $1<p\le q<+\infty$; 
{Following the terminology of \cite{HH,HastoOk}, we give the following definitions. The first three ones concern with the regularity of $\varphi$ with respect to the $t$- variable,  (A1) imposes a bound on how much $\varphi$  can change between
nearby points, while the last one is a continuity assumption with respect to the spatial variable $x$. 
\begin{definition}\label{eleuteri}
Let $p,q>0$. A function $\varphi:\Omega\times[0,+\infty)\to[0,+\infty)$ satisfies
\begin{itemize}
\item[\normalfont(inc)$_p$]\label{inc} 
 if $t\in(0,+\infty)\mapsto\frac{\varphi(x,t)}{t^p}$ is increasing for every $x\in\Omega$
\item[\normalfont(dec)$_q$]\label{dec} if $t\in (0,+\infty)\mapsto\frac{\varphi(x,t)}{t^q}$ is decreasing for every $x\in\Omega$
\item[\normalfont(A0)] \label{azero} if there exists $L\ge 1$ such that $\frac{1}{L}\le \varphi (x,1)\le L$ for every $x\in\Omega$
\item[\normalfont(A1)] \label{a1} if there exists $L\ge 1$ such that, for any ball $B_r(x_0)\subset\Omega$,
\[
\varphi^+_{r,x_0}(t)\le L\,\varphi^-_{r,x_0}(t), \quad\forall t>0 \,\, \mbox{ such that } \,\,\varphi^-_{r,x_0}(t)\in \left [1, \frac{1}{\mathcal{L}^d(B_r(x_0))} \right].
\]
\item[\normalfont(VA1)] \label{va1} if there exists an increasing continuous function $\omega:[0,+\infty)\to [0,1]$ with $\omega(0)=0$ such that, for any ball $B_r(x_0)\subset\Omega$,
\[
\varphi^+_{r,x_0}(t)\le (1+\omega(r))\varphi^-_{r,x_0}(t), \quad\forall t>0 \,\, \mbox{ such that } \,\,\varphi^-_{r,x_0}(t)\in \left [\omega(r), \frac{1}{\mathcal{L}^d(B_r(x_0))} \right].
\]
\end{itemize}
\end{definition} }

\begin{remark}\label{rem:equivalence}
Note  that assumption \vauno{} implies \auno{}, see \cite[Remark 4.2]{HastoOk}.   By \inc{p}, condition \auno{} implies
\item[\normalfont(A1')] \label{a1} there exists $\beta \in (0,1)$ such that, for any ball $B_r(x_0)\subset\Omega$,
\[
\varphi^+_{r,x_0}(\beta t)\le \varphi^-_{r,x_0}(t), \quad\forall t>0 \,\, \mbox{ such that } \,\,\varphi^-_{r,x_0}(t)\in \left [1, \frac{1}{\mathcal{L}^d(B_r(x_0))} \right]\,,
\]
which in its turn implies \auno{} if \dec{q}, hence a doubling condition, holds. In this setting, conditions \azero{} and \auno{} are invariant under a notion of function equivalence (see \cite[Lemma 4.1.3]{HH}), provided the constant $L$ (or $\beta$) is suitably rescaled. In particular, if $\varphi$ satisfies \inc{p}, \dec{q}, \azero{}, and \auno{}, so does $c\varphi$ for every $c \in \mathbb{R}$.
\end{remark}

{\begin{remark}
If $\varphi$ satisfies \inc{p} (resp., \dec{q}) for some $p>0$ (resp., $q>0$), then so do $\varphi^+_{r,x_0}$ and $\varphi^-_{r,x_0}$ for any $B_{r}(x_0)\subset\Omega$.
\label{rem:scaledvarphi}
\end{remark}}

{\begin{remark} \label{rem:varinc}
If $\varphi:\Omega\times [0,+\infty)\to[0,+\infty)$ is convex and $\varphi(x,0)=0$ for every $x\in\Omega$, then $\varphi$ satisfies \inc{1}. If $\varphi$ satisfies \inc{p_1}, then it satisfies \inc{p_2} for every $0<p_2\leq p_1$. If $\varphi$ satisfies \dec{q_1}, then it satisfies \dec{q_2} for every $q_2\geq q_1$.
\end{remark}}

Next simple results {can be found in \cite[Section 3]{HastoOk}.}

{\begin{proposition}\label{prop:properties}
Let $1<p\le q<+\infty$ and $\varphi\in \Phic(\Omega)$ with right derivative $\varphi_t$. Assume that $\varphi_t$ satisfies \inc{p-1} and \dec{q-1}. Then
\begin{enumerate}
\item[$(i)$] $\varphi$ satisfies \inc{p} and \dec{q}, and the following estimate hold:
\begin{equation}\label{cons2}
\varphi(x,s)\min\{t^p,t^q\}\le\varphi(x,ts)\le\max\{t^p,t^q\}\varphi(x,s),\quad \forall x\in\Omega, \,\, \forall s,t\in[0,+\infty).
\end{equation}
\item[(ii)] $\varphi(x,t)$ and $t\varphi_t(x,t)$ are equivalent, in the sense that \begin{equation}\label{cons1}
p\,\varphi(x,t)\le t\,\varphi_t(x,t)\le q\,\varphi(x,t),\quad \forall (x,t)\in\Omega\times[0,+\infty);
\end{equation}
\item[(iii)] if, in addition, $\varphi_t$ complies with \azero{}, then also $\varphi$ does with constants depending on $L,p,q$. More precisely, 
\begin{equation}\label{v0phi}
\frac{1}{Lq}\le\varphi(x,1)\le \frac{L}{p}, \ \ \ \forall x\in\Omega.
\end{equation}
\end{enumerate}
If, in addition, $\varphi(x, \cdot )\in C^1([0,+\infty))$ for every $x\in\Omega$, then $\varphi\in N(\Omega)$.
\end{proposition} 
}

{
For $\varphi\in\Phiw(\Omega)$, the \textit{generalized Orlicz space} is defined by 
\[
L^{\varphi}(\Omega):=\big\{f\in L^0(\Omega):\|f\|_{L^\varphi(\Omega)}<\infty\big\}
\] 
with the (Luxemburg) norm 
\[
\|f\|_{L^\varphi(\Omega)}:=\inf\bigg\{\lambda >0: \varrho_{\varphi}\Big(\frac{f}{\lambda}\Big)\leq 1\bigg\},
\ \ \text{where}\ \ \varrho_{\varphi}(f):=\int_\Omega\varphi(x,|f(x)|)\,\mathrm{d}x.
\]
We denote by $W^{1,\varphi}(\Omega)$ the set of $f\in L^{\varphi}(\Omega)$ satisfying that $\partial_1f,\dots,\partial_df \in L^{\varphi}(\Omega)$, where $\partial_if$ is the weak derivative of $f$ in the $x_i$-direction, with the norm $\|f\|_{W^{1,\varphi}(\Omega)}:=\|f\|_{L^\varphi(\Omega)}+\sum_i\|\partial_if\|_{L^\varphi(\Omega)}$. Note that if $\varphi$ satisfies \dec{q} for some $q\ge 1$, then $f\in L^\varphi(\Omega)$ if and only if $\varrho_\varphi(f)<\infty$, and if $\varphi$ satisfies \azero{}, \inc{p} and \dec{q} for some $1<p\leq q$, then $L^\varphi(\Omega)$ and $W^{1,\varphi}(\Omega)$ are reflexive Banach spaces. In addition we denote by $W^{1,\varphi}_0(\Omega)$ the closure of $C^\infty_0(\Omega)$ in $W^{1,\varphi}(\Omega)$.}

The following version of Sobolev-Poincarè inequality for weak $\Phi$-functions can be deduced by \cite[Proposition 3.6]{HHK}.

\begin{proposition}
Let $B_r\subset\mathbb{R}^d$ be a ball and $\varphi\in\Phiw(B_r)$ be complying with \azero{}, \auno{}, \inc{p}, \dec{q}, $1\leq p<q$, and let $s\in[1,p]$ with $s<\frac{d}{d-1}$. Then there exists a constant $C_{\rm P}=C_{\rm P}(d,s,p,q,L)$ such that
\begin{equation}
\dashint_{B_r} \varphi\left(x, \frac{|u-(u)_{B_r}|}{2r}\right)\,\mathrm{d}x \leq C_{\rm P} \left( \left ( \dashint_{B_r}\varphi(x,|\nabla u|)^\frac{1}{s}\,\mathrm{d}x \right)^s + 1 \right)
\label{eq:poincare1}
\end{equation}
for any $u\in W^{1,1}(B_r)$ such that $\|\nabla u\|_{L^\varphi(B_{r})}\leq1$. 

%{If, in addition, $\varphi$ satisfies \dec{q} and %if $u\in W^{1,1}(B_r)$ such that $\|\nabla u\|_{L^\varphi(B_{r})}\leq1$ and $u=0$ on $\partial B_r$, or 
%$u=0$ on $A\subset B_r$ with $\mathcal{L}^d(A)>0$, then
%\begin{equation}
%\dashint_{B_r} \varphi\left(x, \frac{|u|}{r}\right)\,\mathrm{d}x \leq 2^{q-1}C_{\rm P}\left(1+\left(\frac{\mathcal{L}^d(B_r)}{\mathcal{L}^d(A)}\right)^q\right) \left( \left ( \dashint_{B_r}\varphi(x,|\nabla u|)^\frac{1}{s}\,\mathrm{d}x \right)^s + 1 \right)\,.
%\label{eq:poincare2}
%\end{equation}}
%where in the second case the constant $c$ depends also on $\frac{\mathcal{L}^d(B_r)}{\mathcal{L}^d(A)}$.  }
\label{prop:sobpoinc}
\end{proposition}
\proof
The role of the further assumption \dec{q} with respect to \cite[Proposition 3.6]{HHK} is that it allows to transfer the constant $\beta_3=\beta_3(d,s,\varphi)>0$ therein to the right-hand side of the inequality as $C_P:=\max\{\beta_3^{-q},1\}$. 
\endproof

\subsection{Regularized Orlicz function} \label{sec:regularizedphi}

Let $\varphi \in \Phic(\Omega)$, $\varphi(x,\cdot)\in C^1([0,\infty))$ satisfying \auno{}, and such that $\varphi_t$ comply with \azero{}, \inc{p-1} and \dec{q-1} for some $1<p\leq q$. Then, as proven in \cite[Proposition 5.10]{HastoOk}, on each ball $B=B_{2r}(x_0)$ a regularized function $\widetilde{\varphi}=\widetilde{\varphi}_B\in C^1([0,\infty))\cap C^2((0,\infty))$ can be constructed such that
\begin{itemize}
\item[(i)] $\widetilde{\varphi}$ satisfies \azero{}, \inc{p}, \dec{q}, while $\widetilde{\varphi}^\prime$ complies with \azero{}, \inc{p-1} and \dec{q-1}. In particular,
\begin{equation}
\widetilde{\varphi}^\prime(t) \sim t \widetilde{\varphi}^{\prime\prime}(t)\,, \quad \mbox{ uniformly for all $t>0$.}
\label{eq:equivtphi}
\end{equation}
\item[(ii)] \begin{equation}
\widetilde{\varphi}(t) \leq c (\varphi(x,t)+1) \quad \mbox{ for all } (x,t)\in B\times [0,\infty)\,.
\label{eq:comptildephi}
\end{equation} 
\end{itemize}

%Correspondingly, we introduce the map $\mathbf{V}: \mathbb{R}^d \to \mathbb{R}^d$ given by
%\begin{equation}
%\label{Vmap}
%\mathbf{V}(z):= \sqrt{\frac{\widetilde{\varphi}^\prime(|z|)}{|z|}}\,  z \,.
%\end{equation}
%The following properties
%\begin{equation}
%|\mathbf{V}(z_1)-\mathbf{V}(z_2)|^2 \sim \widetilde{\varphi}^{\prime\prime}(|z_1|+|z_2|) |z_1-z_2|^2\,, \quad \mbox{ for every $z_1,z_2\in\mathbb{R}^d$}
%\label{eq:equivtphi2}
%\end{equation}
%\begin{equation}
%|\mathbf{V}(z)|^2 \sim \widetilde{\varphi}(|z|)\,, \quad \mbox{ for every $z\in\mathbb{R}^d$,}
%\label{eq:equivtphi4}
%\end{equation}
%hold with uniform constants with respect to $z,z_1,z_2$. They are nowadays standard and can be found in \cite[Lemma 3]{DIEETT}. 
Moreover, taking into account \eqref{eq:equivtphi}, it can be shown that (see \cite[Lemma 3.8(2)-(3)]{HastoOk})
\begin{equation}
\widetilde{\varphi}^{\prime\prime}(|z_1|+|z_2|) |z_1-z_2|^2 \lesssim \widetilde{\varphi}(|z_1|)-\widetilde{\varphi}(|z_2|)  - \frac{\widetilde{\varphi}^\prime(|z_2|)}{|z_2|} z_2 \cdot (z_1-z_2)\,, \quad \mbox{ for every $z_1,z_2\in\mathbb{R}^d$. }
\label{eq:equivtphi3}
\end{equation}
and
\begin{equation}
\widetilde{\varphi}(|z_1-z_2|) \lesssim 
\varepsilon\left[\widetilde{\varphi}(|z_1|)+\widetilde{\varphi}(|z_2|)\right]+ \varepsilon^{-1}\widetilde{\varphi}''(|z_1|+|z_2|)|z_1-z_2|^2\,, \quad \mbox{ for every $z_1,z_2\in\mathbb{R}^d$, $\varepsilon>0$. }
\label{eq:equivtphi5}
\end{equation}
%We also recall from \cite[Lemma~A.2]{DiKaSch} that for any ball $B\subset\Omega$ and $f \in W^{1,\widetilde{\varphi}}(B)$,
%\begin{equation}
%\dashint_{B}|{\bf V}(f)-{\bf V}((f)_{B})|^2\,\mathrm{d}x \sim  \dashint_{B}|{\bf V}(f)-({\bf V}(f))_{B}|^2\,\mathrm{d}x \sim \dashint_{B} \widetilde{\varphi}_{|(f)_{B}|}(|f-(f)_{B}|)\,\mathrm{d}x\,.
%\label{eq:equivalencebis}
%\end{equation} 

\subsection{Some technical lemmas} \label{sec:techlemmas}
The following lemma, useful in order to re-absorb certain terms, is a variant of the classical \cite[Lemma~6.1]{GIUSTI}.

\begin{lemma}{(\cite[Lemma 4.3]{HHL})}\label{lem:iteration}
Let $Z$ be a bounded non-negative function in the interval $[r,R]$ and let $X$ be an almost decreasing function on $[0,+\infty)$. Assume that there exists $\theta\in[0,1)$ such that
\begin{equation*}
Z(s)\leq\theta Z(t) + X\left(\frac{1}{t-s}\right)\,,
\end{equation*}
for all $r\leq s < t \leq R$. Then 
\begin{equation*}
Z(r) \lesssim X\left(\frac{1}{R-r}\right)\,,
\end{equation*}
where the implicit constant depends on the constant of almost monotonicity and on $\theta$. 
\end{lemma}

In order to derive reverse H\"older estimates, we need a variant of the results by Gehring~\cite{Gehring} and Giaquinta-Modica~\cite[Theorem~6.6]{GIUSTI}.
\begin{lemma}\label{lem:gehring}
Let $B_0\subset\mathbb{R}^n$ be a ball, $f\in L^1(B_0)$, and $g\in L^{\sigma_0}(B_0)$ for some $\sigma_0>1$. Assume that for some $\theta\in(0,1)$, $c_1>0$ and all balls $B$ with $2B\subset B_0$
\begin{equation*}
\dashint_B |f|\,\mathrm{d}x\leq c_1 \left(\dashint_{2B}|f|^\theta\,\mathrm{d}x\right)^{1/\theta} + \dashint_{2B}|g|\,\mathrm{d}x\,.
\end{equation*}
Then there exist $\sigma_1>1$ and $c_2>1$ such that 
%$g\in L^{\sigma_1}_{\rm loc}(B)$ and 
for all $\sigma_2\in[1,\sigma_1]$
\begin{equation*}
\left(\dashint_{B}|f|^{\sigma_2}\,\mathrm{d}x\right)^{1/{\sigma_2}}\leq c_2 \dashint_{2B}|f|\,\mathrm{d}x + c_2 \left(\dashint_{2B}|g|^{\sigma_2}\,\mathrm{d}x\right)^{1/{\sigma_2}}\,.
\end{equation*}
\end{lemma}

The following iteration lemma can be found, e.g., in \cite[Lemma 7.3]{GIUSTI}. %\cite[Lemma 7.1]{HastoOk} (see also \cite[Lemma 2.1 in Chapter III]{Giaquinta}). 
\begin{lemma}
Let $f:[0,R] \to [0,\infty)$ be a non-decreasing function. Assume that 
$$
f(\rho)\leq A\left(\left(\frac{\rho}{r}\right)^\delta+\varepsilon\right)f(r)+B r^\gamma \quad \mbox{ for all $0<\rho\leq r\leq R$,}
$$ 
for positive constants $A$ and $B$, and $\delta\leq \gamma$. Then for any $\sigma \in(0,\delta)$, there exist $\varepsilon_0, c>0$ depending only on $\gamma$, $\delta$, $A$ and $\sigma$ such that if $\varepsilon<\varepsilon_0$, then
$$
f(\rho)\leq c\left(\left(\frac{\rho}{r}\right)^{\delta-\sigma} f(r)+B\rho^{\delta-\sigma}\right)\,.
$$ 
\label{lem:iterationlemma}
\end{lemma}

\section{Local Lipschitz regularity of almost minimizers} \label{sec:intregularity}

\subsection{Regularity estimates for autonomous problems} %: the $\widetilde{\varphi}$-harmonic replacement}
\label{sec:regauto}

{Let $\psi\in\Phic\cap C^1([0,\infty))\cap C^2((0,\infty))$ with $\psi^\prime$ satisfying \inc{p-1} and \dec{q-1} for some $1<p\leq q$.  %be the $N$-function defined in Section \ref{sec:regularizedphi}. 
For a given ball $B_r(x_0)\Subset \Omega$ and $w_0\in W^{1,\psi}(B_r(x_0))$ %local almost-minimizer of $\mathcal{F}$, 
%and a , 
we consider a weak solution to the Dirichlet problem
\begin{equation}
\begin{cases}
\displaystyle {\rm div} \left (\frac{\psi^\prime(|\nabla w|)}{|\nabla w|}\nabla w \right) = 0 & \mbox{ \,\, in \,\, $ B_r(x_0)$,  }  \\
w = w_0  & \mbox{ \,\, on \,\, $\partial B_r(x_0)$.  } 
\end{cases}
\label{eq:systemphi}
\end{equation}
%or, equivalently, a solution %$u^*_r\in u+ W^{1,\widetilde{\varphi}}_0(B_r(x_0))$ 
%to the minimization problem
%\begin{equation}
%\min_{w\in w_0+ W^{1,\psi}_0(B_r(x_0))} \int_{B_r(x_0)}\psi(|\nabla w|)\,\mathrm{d}x \,.
%\label{eq:functionalphi}
%\end{equation}

As proven in \cite[Lemma 4.12]{HastoOk},  %\cite[Lemma 5.8 and Lemma 6.4]{DSV}, 
the following Harnack-type inequality and excess decay estimate hold for any such $w$.
\begin{proposition}
Let $\psi\in\Phic\cap C^1([0,\infty))\cap C^2((0,\infty))$ with $\psi^\prime$ satisfying \inc{p-1} and \dec{q-1} for some $1<p\leq q$. Let $w\in W^{1,\psi}(B_r(x_0))$ be a weak solution to \eqref{eq:systemphi}. Then there exists $\mu_0=\mu_0(d,p,q)\in(0,1)$ such that $\nabla w \in C^{0,\mu_0}_{{\rm loc}}(B_r(x_0);\mathbb{R}^d)$ and the following estimates hold: there exists a constant $c=c(d,p,q)>0$ such that, for every $B_\rho(y)\subset B_r(x_0)$,
\begin{equation}
\sup_{B_{\rho/2}(y)} |\nabla w| \leq c \dashint_{B_\rho(y)} |\nabla w|\,\mathrm{d}x \,,
\label{eq:5.15DSV}
\end{equation}
and for any $\tau\in(0,1)$,
\begin{equation}
\dashint_{B_{\tau \rho}(y)} \left| \nabla w - (\nabla w)_{B_{\tau \rho}(y)} \right| \, \mathrm{d}x \leq c \tau^{\mu_0} \dashint_{B_{\rho}(y)} \left| \nabla w \right|\, \mathrm{d}x \,.
\label{eq:DSVestim}
\end{equation}
\label{lem:lem4.12HOK}
\end{proposition}

For a given $u\in W^{1,\varphi}_{\rm loc}(\Omega)$ local almost-minimizer of $\mathcal{F}$, and a ball $B_r(x_0)\Subset \Omega$, we consider $\widetilde{\varphi}$ the $N$-function defined in Section \ref{sec:regularizedphi} on $B_r(x_0)$, and the unique weak solution to the Dirichlet problem
\begin{equation}
\begin{cases}
\displaystyle {\rm div} \left (\frac{\widetilde{\varphi}^\prime(|\nabla w|)}{|\nabla w|}\nabla w \right) = 0  & \mbox{ \,\, in \,\, $ B_r(x_0)$,  } \\
w = u & \mbox{ \,\, on \,\, $\partial B_r(x_0)$, } 
\end{cases}
\label{eq:systemphitilde}
\end{equation}
or, equivalently, the solution %$u^*_r\in u+ W^{1,\widetilde{\varphi}}_0(B_r(x_0))$ 
to the minimization problem
\begin{equation}
\min_{w\in u+ W^{1,\widetilde{\varphi}}_0(B_r(x_0))} \int_{B_r(x_0)}\widetilde{\varphi}(|\nabla w|)\,\mathrm{d}x \,.
\label{eq:functionalphi}
\end{equation}

The existence and uniqueness in the minimization problem above follows from the fact that, by \eqref{eq:comptildephi}, $u\in W^{1,\widetilde{\varphi}}(B_r(x_0))$, for which $u$ is an admissible boundary-value function.  This suggests the following definition. 
\begin{definition}
We define the \emph{$\widetilde{\varphi}$-harmonic replacement} of $u$ in $B_r(x_0)$, and we denote it by $v_r$, as the unique solution to the variational problem  \eqref{eq:functionalphi}.
\end{definition}
Since we may choose $\psi=\widetilde{\varphi}$ in Proposition \ref{lem:lem4.12HOK}, we have $\nabla v_r \in C^{0,\mu_0}_{\rm loc}(B_r(x_0);\mathbb{R}^d)$ and the estimates \eqref{eq:5.15DSV} and \eqref{eq:DSVestim} hold for $v_r$ on every $B_\rho(y)\subset B_r(x_0)$.

\subsection{Preliminary regularity results for almost-minimizers and comparison estimates}\label{sec:compar}%Caccioppoli inequality, higher integrability and comparison estimates}

We start by proving a Caccioppoli-type inequality for almost minimizers of $\mathcal{F}$. 

\begin{lemma} %Let $\varphi\in \Phiw(\Omega)$ be doubling. 
Let $\varphi\in \Phic(\Omega)$ be such that \dec{q} holds for some $q>0$. Let $u$ be an almost-minimizer of $\mathcal{F}$ in $\Omega$, with constant $\kappa\leq \kappa_0$ and exponent $\beta$, and let $x_0 \in \Omega$ and $B_{2r}(x_0) \Subset \Omega$, with $2r\leq 1$. Then there exists a constant $c=c(q,\kappa_0)$ such that
\begin{equation}
\dashint_{B_{r}(x_0)} \varphi(x,|\nabla u|)\,\mathrm{d}x  \leq c \left( \dashint_{B_{2r}(x_0)} \varphi\left(x,\frac{|u-(u)_{x_0,2r}|}{2r}\right)\,\mathrm{d}x + \lambda \right) \,.
\end{equation}
\label{lem:caccioppoli}
\end{lemma}

\begin{proof}
Let $1\leq s < t \leq 2$ and $\eta\in C^\infty_0(B_{tr}(x_0),[0,1])$ be a cut-off function such that $\eta\equiv 1 $ on $B_{sr}(x_0)$ and $|\nabla \eta|\leq \frac{2}{(t-s)r}$. Set $v:=u-\eta(u-(u)_{x_0,2r})$. Then we have $v=u$ on $\partial B_{tr}(x_0)$ and
\begin{equation*}
\nabla v = (1-\eta) \nabla u - \nabla \eta ( u-(u)_{x_0,2r})\,. 
\end{equation*}
Then, since $u$ is an almost minimizer of $\mathcal{F}$, using also the convexity of $\varphi(x,\cdot)$ and \dec{q} we get
\begin{equation*}
\begin{split}
\int_{B_{sr}(x_0)} \varphi(x,|\nabla u|)\,\mathrm{d}x & \leq  \mathcal{F}(u; B_{tr}(x_0)) \\
& \leq \left(1+\kappa_0 \right)  \mathcal{F}(v; B_{tr}(x_0)) \\
& \leq c_1\left(\int_{B_{tr}(x_0)} \varphi(x,|\nabla u|(1-\eta))\,\mathrm{d}x + \int_{B_{tr}(x_0)} \varphi\left(x,\frac{|u-(u)_{x_0,2r}|}{(t-s)r}\right)\,\mathrm{d}x + \lambda\mathcal{L}^d(B_{2r})\right) \\
& \leq c_1 \left(\int_{B_{tr}(x_0)\backslash B_{sr}(x_0)} \varphi(x,|\nabla u|)\,\mathrm{d}x +  \int_{B_{2r}(x_0)} \varphi\left(x,\frac{|u-(u)_{x_0,2r}|}{(t-s)r}\right)\,\mathrm{d}x + \lambda\mathcal{L}^d(B_{2r})\right) \,,
\end{split}
\end{equation*}
where $c_1=c_1(q,\kappa_0)$, whence adding $c_1\int_{B_{sr}(x_0)} \varphi(x,|\nabla u|)\,\mathrm{d}x$ to both the sides and then dividing by $1+c_1$ we get
\begin{equation*}
\begin{split}
\int_{B_{sr}(x_0)} \varphi(x,|\nabla u|)\,\mathrm{d}x \leq \frac{c_1}{1+ c_1} \left(\int_{B_{tr}(x_0)} \varphi(x,|\nabla u|)\,\mathrm{d}x +\int_{B_{2r}(x_0)} \varphi\left(x,\frac{|u-(u)_{x_0,2r}|}{(t-s)r}\right)\,\mathrm{d}x + \lambda \mathcal{L}^d(B_{2r})\right) \,.
\end{split}
\end{equation*}
Now, 
%a standard iteration procedure based on 
an application of Lemma \ref{lem:iteration} with $\theta:=\frac{c_1}{1+ c_1}$, $Z(t):=\int_{B_{tr}(x_0)} \varphi(x,|\nabla u|)\,\mathrm{d}x$ and $X(\tau):= \frac{c_1}{1+c_1} \int_{B_{2r}(x_0)} \left(\varphi\left(x, \tau |u-(u)_{x_0,2r}|\right) +\lambda\right)\,\mathrm{d}x$ gives 
\begin{equation*}
\begin{split}
\int_{B_{r}(x_0)} \varphi(x,|\nabla u|)\,\mathrm{d}x \leq c\int_{B_{2r}(x_0)} \varphi\left(x,\frac{|u-(u)_{x_0,2r}|}{2r}\right)\,\mathrm{d}x + c \lambda \mathcal{L}^d(B_{2r}) \,.
\end{split}
\end{equation*}
Taking the average of both the sides we finally get the desired result. 
\endproof

\end{proof}

The following lemma contains a higher integrability result and reverse H\"{o}lder type estimates for the gradient of an almost minimizer of $\mathcal{F}$. 

\begin{lemma} Let $\varphi\in\Phiw(\Omega)$ satisfy \azero{}, \auno{}, \inc{p}, \dec{q} with constant $L\geq1$ and $1<p\leq q$. Let $u\in W^{1,\varphi}_{\rm loc}(\Omega)$ be an almost-minimizer of $\mathcal{F}$ in $\Omega$ with constant $\kappa\leq\kappa_0$ and exponent $\beta$, and let $x_0 \in \Omega$ and $B_{2r}(x_0) \Subset \Omega$, with $\|\nabla u\|_{L^\varphi(B_{2r}(x_0))}\leq1$, and $2r\leq 1$. Then 
\begin{enumerate}
\item[(i)](Higher integrability) there exist $s_0=s_0(d,p,q,L)>0$ and $c=c(d, p, q, L,\kappa_0)\geq1$ such that
\begin{equation}
\left(\dashint_{B_{r}(x_0)} \varphi(x,|\nabla u|)^{1+s_0}\,\mathrm{d}x \right)^{\frac{1}{1+s_0}} \leq c 2^{\frac{ds_0}{1+s_0}}\delta^{-\frac{ds_0}{1+s_0}} \left(\dashint_{B_{(1+\delta)r}(x_0)} \varphi\left(x,|\nabla u|\right)\,\mathrm{d}x +\Lambda \right) \,,
\end{equation}
for any $\delta\in(0,1]$, where $\Lambda:=\lambda+1$. In particular, this implies $\varphi(\cdot,|\nabla u|) \in L^{1+s_0}_{\mathrm{loc}}(\Omega)$.
\item[(ii)](Reverse H\"older type estimates) for every $t\in(0,1]$, there exist $c_t=c_t(d, p, q, L,\kappa_0, t)>0$ such that
\begin{equation}
\left(\dashint_{B_{r}(x_0)} \varphi(x,|\nabla u|)^{1+s_0}\,\mathrm{d}x \right)^{\frac{1}{1+s_0}} \leq c_t \left(\left(\dashint_{B_{2r}(x_0)} \varphi\left(x,|\nabla u|\right)^t\,\mathrm{d}x\right)^\frac{1}{t} +\Lambda \right) \,,
\label{eq:reverseholderr}
\end{equation}
and $c=c(d, p, q, L,\kappa_0)\geq1$ such that
\begin{equation}
\dashint_{B_{r}(x_0)} \varphi(x,|\nabla u|)\,\mathrm{d}x \leq \left(\dashint_{B_{r}(x_0)} \varphi(x,|\nabla u|)^{1+s_0}\,\mathrm{d}x \right)^{\frac{1}{1+s_0}} \leq c \left(\varphi^-_{B_{2r}(x_0)}\left(\dashint_{B_{2r}(x_0)} |\nabla u|\,\mathrm{d}x\right) +\Lambda \right) \,.
\end{equation}
\end{enumerate} 
\label{lem:highint}
\end{lemma}
\begin{proof}
The higher integrability result for $\nabla u$ in $(i)$ can be obtained in a standard way combining the Caccioppoli inequality of Lemma \ref{lem:caccioppoli} with the Sobolev-Poincarè inequality (Proposition \ref{prop:sobpoinc}), by Gehring's lemma (Lemma \ref{lem:gehring}). The reverse H\"{o}lder type inequalities $(ii)$ follow from $(i)$ by a similar argument as for \cite[Lemma 4.7]{HastoOk}. We omit further details. 
\end{proof}

Note that, under our assumption on $\varphi$, it holds that
\begin{equation*}
\|\nabla u\|_{L^\varphi(B_{2r}(x_0))}\leq1 \iff \int_{B_{2r}(x_0)} \varphi(x,|\nabla u|)\,\mathrm{d}x \leq 1
\end{equation*}
(it is sufficient that $\varphi\in\Phiw(\Omega)$ and $\varphi(x,\cdot)$ be left-continuous; see, e.g., \cite[Lemma 3.2.3]{HH}). 

\begin{remark}[Choice of small radii] \label{rmk:smallness}
(i) Arguing as in \cite{HastoOk}, with fixed $\Omega^\prime \Subset\Omega$ we can show that there exists $r_0 \in (0,1)$, $r_0=r_0(d,L,\omega(\cdot),\|\varphi(\cdot, |\nabla u|)^{1+s_0}\|_{L^1(\Omega')})$ satisfying
\begin{equation}
r_0\leq \frac{1}{2}\,, \quad \omega(2r_0)\leq \frac{1}{L}\,, \quad \mathcal{L}^d(B_{2r_0}) \leq \min\left \{\frac{1}{2L}, 2^{-\frac{2(1+s_0)}{s_0}} \left(\int_{\Omega^\prime} \varphi(x, |\nabla u|)^{1+s_0}\,\mathrm{d}x\right)^{-\frac{2+s_0}{s_0}} \right\}\,,
\label{eq:smallradius}
\end{equation}
where $L\geq1$ is that of condition \azero{} and $s_0$ is the exponent of Lemma \ref{lem:highint}, 
%Lemma \ref{lem:highint} implies $\varphi(\cdot,|\nabla u|) \in L^{1+s_0}_{\mathrm{loc}}(\Omega)$. 
%Then, arguing as in \cite[(5.3)-(5.4)]{HastoOk}, as a consequence of \eqref{eq:smallradius}, H\"{older}'s inequality and Young's inequality,  for each $\Omega'\Subset \Omega$, there exists $r_0 \in (0,\frac{1}{2}]$, $r_0=r_0(d,L,\omega(\cdot),\|\varphi(\cdot, |\nabla u|)^{1+s_0}\|_{L^1(\Omega')})$ 
such that  for any $B_{2r}
(x_0)\subset \Omega'$ with $r\in(0,r_0]$, it holds that
\begin{equation*}\label{eq:assless1}
\int_{B_{2r}(x_0)} \varphi(x,|\nabla u|)\,\mathrm{d}x \le 1 \,.
\end{equation*}
Therefore, we can exploit the estimates of Lemma \ref{lem:highint} in each of these balls.  \\
\noindent
(ii) If, in addition, $u\in L^\infty_{\rm loc}(\Omega)$, the choice of $s_0$ and, accordingly, of $r_0$ can be done in such a way that the dependence of $r_0$ on $u$ is through $\|u\|_{L^\infty(\Omega')}$.  This observation will be crucial when dealing with some auxiliary sequences in our proofs.
To show this, let $\Omega''\Subset\Omega'\Subset\Omega$ be fixed, $y\in\overline{\Omega''}$ and $r< \frac{1}{4}\min\{1,{\rm dist}(\Omega'',\partial\Omega')\}$. Then, by the Caccioppoli inequality of Lemma \ref{lem:caccioppoli}, \inc{p}, \dec{q}, \azero \  we first find
\begin{equation*}
\begin{split}
\int_{B_{2r}(y)} \varphi(x,|\nabla u|)\,\mathrm{d}x & \leq c(q,\kappa_0)\left(\int_{B_{4r}(y)} \varphi^+_{B_{4r}(y)}\left(\frac{\|u\|_{L^\infty(B_{4r})}}{2r}\right)\,\mathrm{d}x + \lambda \mathcal{L}^d(B_{r})\right) \\
& \leq c(d,p,q,\kappa_0,L) r^d\left(\max\left\{\left({\textstyle\frac{\|u\|_{L^\infty(\Omega')}}{r}}\right)^p, \left({\textstyle\frac{\|u\|_{L^\infty(\Omega')}}{r}}\right)^q\right\}+\lambda\right) \\
& =: \bar{c}(d,p,q,L,\kappa_0,r,\|u\|_{L^\infty(\Omega')},\lambda) \,,
\end{split}
\end{equation*}
where $\bar{c}$ depends increasingly on $\|u\|_{L^\infty(\Omega')}$.
If $\bar{c}\leq 1$, then by \eqref{eq:reverseholderr} for $t=1$ and a covering argument we get the bound
\begin{equation}
\begin{split}
\int_{\Omega''} \varphi(x,|\nabla u|)^{1+s_0}\,\mathrm{d}x  & \leq [c_1 \left(\bar{c} +\Lambda \right)]^{1+s_0} \\
& =: \tilde{c} = \tilde{c}(d,p,q,L,\kappa_0,\|u\|_{L^\infty(\Omega')},{\rm dist}(\Omega'',\partial\Omega'),\lambda)\,,
\end{split}
\label{eq:newargument}
\end{equation}
where $\tilde{c}$ depends increasingly on $\|u\|_{L^\infty(\Omega')}$. Hence, in \eqref{eq:smallradius}, it will suffice to require 
\begin{equation*}
r\leq \frac{1}{2}\,, \quad \omega(2r)\leq \frac{1}{L}\,, \quad \mathcal{L}^d(B_{2r}) \leq \min\left \{\frac{1}{2L}, 2^{-\frac{2(1+s_0)}{s_0}} \tilde{c}^{-\frac{2+s_0}{s_0}} \right\}\,.
\end{equation*}
If $\bar{c}>1$, we can apply the same argument to $\bar{\varphi}(x,t):=\frac{1}{\bar{c}}\varphi(x,t)$, observing that \azero, \inc{p}, \dec{q} and \auno{} still hold (see Remark \ref{rem:equivalence}).  In this case, as the constants in \azero{} may have changed, the value of $s_0>0$ could become possibly smaller.
\end{remark}

With the following Lemma, we prove a further reverse H\"{o}lder type inequality for $\nabla u$ and a Calderón-Zygmund type estimate for the problem \eqref{eq:systemphitilde}. 

\begin{lemma}\label{lem:lemma5.15HOK}
Let $\varphi$ be as in Lemma \ref{lem:highint}. Let $u\in W^{1,\varphi}_{\rm loc}(\Omega)$ be an almost minimizer of $\mathcal{F}$ in $\Omega$ with constant $\kappa\leq\kappa_0$ and exponent $\beta$, and $v_r\in W^{1,\widetilde{\varphi}}(B_r(x_0))$ be the $\widetilde{\varphi}$-harmonic replacement of $u$ in $B_r(x_0)$, where $B_{2r}(x_0)\Subset \Omega$ with $r$ satisfying \eqref{eq:smallradius}, and $\widetilde{\varphi}$ defined on $B_r(x_0)$ as in Section \ref{sec:regularizedphi}. Then there exists a constant $c = c(d, p, q, L,\kappa_0)\geq1$ such that
\begin{equation}\label{eq:5.16HOK}
\left(\dashint_{B_r(x_0)}\varphi(x,|\nabla u|)^{1+s_0}\,\mathrm{d}x\right)^{\frac{1}{1+s_0}} 
\le c\left[
\widetilde{\varphi} \left(\dashint_{B_{2r}(x_0)}|\nabla u|\,\mathrm{d}x\right)+\Lambda\right]
\end{equation}
and
\begin{equation}\label{eq:5.17HOK}\begin{aligned}
\dashint_{B_r(x_0)}\varphi(x,|\nabla v_r|)\,\mathrm{d}x 
&\leq \left(\dashint_{B_r(x_0)}\varphi(x,|\nabla v_r|)^{1+\frac{s_0}{2}}\,\mathrm{d}x\right)^{\frac{1}{1+\frac{s_0}{2}}}\\
&\le 
c \left(\dashint_{B_r(x_0)}\varphi(x,|\nabla u|)^{1+\frac{s_0}{2}}\,\mathrm{d}x+\Lambda\right)^{\frac{1}{1+\frac{s_0}{2}}}
\end{aligned}\end{equation}
Moreover,
\begin{equation}\label{eq:5.18HOK}
\dashint_{B_r(x_0)}| \nabla v_r |\,\mathrm{d}x 
\le c\left(
\dashint_{B_{2r}(x_0)}|\nabla u|\,\mathrm{d}x+\Lambda\right).
\end{equation}
\end{lemma}

\begin{proof}
The proof can be obtained exactly as in \cite[Lemma 5.15]{HastoOk} by using \azero{}, Lemma \ref{lem:highint}$(ii)$ and some general techical results (see \cite[Proposition 5.12 and Lemma 4.15]{HastoOk}). We omit further details.
\end{proof}

With the following Proposition, we establish a comparison estimate between the gradient of an almost minimizer of $\mathcal{F}$ and that of its $\widetilde{\varphi}$-harmonic replacement in a ball. From now on, the stronger condition \vauno{} is needed in place of \auno{}.

\begin{proposition}\label{prop:stimav}
Let $\varphi \in \Phic(\Omega)$, $\varphi(x,\cdot)\in C^1([0,\infty))$ be satisfying \vauno{}, such that $\varphi_t$ satisfies \azero{}, \inc{p-1} and \dec{q-1} for some $1<p\leq q$. Let $u\in W^{1,\varphi}_{\rm loc}(\Omega)$ be an almost-minimizer of $\mathcal{F}$ in $\Omega$, with constant $\kappa\leq \kappa_0$ and exponent $\beta$, let $x_0 \in \Omega$ and $r>0$ be such that $B_{2r}(x_0) \Subset \Omega$ and complying with \eqref{eq:smallradius}. Let $\widetilde{\varphi}$ be defined on $B_r(x_0)$ as in Section \ref{sec:regularizedphi}, 
and $v_r\in W^{1,\widetilde{\varphi}}(B_r(x_0))$ be the $\widetilde{\varphi}$-harmonic replacement of $u$ in $B_r(x_0)$. Then, there exists a constant $c=c(d,p,q,L, \kappa_0)\geq1$ such that
\begin{equation}
\dashint_{B_r(x_0)} \widetilde{\varphi}^{\prime\prime}(|\nabla u|+|\nabla v_r|) |\nabla u - \nabla v_r|^2\, \mathrm{d}x \lesssim \left( \omega(2r)^\frac{p}{q} + r^{\min\{\beta,\gamma\}} \right)  \widetilde{\varphi}\left(\dashint_{B_{2r}(x_0)} \left|\nabla u \right| \, \mathrm{d}x\right) + \Lambda  \,. 
\label{eq:stimaV}
\end{equation}
%\begin{equation}
%\dashint_{B_r(x_0)} \left| \nabla u - \nabla v_r \right| \, \mathrm{d}x \leq c \left(\omega(2r)^\frac{p}{2q^2} + r^\frac{\min\{\beta,\gamma\}}{2q} \right) \left(\dashint_{B_{2r}(x_0)} \left|\nabla u \right| \, \mathrm{d}x + \Lambda \right) \,.
%\label{eq:stimasenzaV}
%\end{equation} 
where $\gamma:=\min\left\{1,\frac{d s_0^2}{4(2+s_0)}\right\}$ and $s_0$ is the exponent of Lemma \ref{lem:highint}. 
\label{lem:stimaV}
\end{proposition}

\begin{proof}
%{\bf Step 1:} We first prove that
%\begin{equation}
%\int_{B_r(x_0)} \widetilde{\varphi}^{\prime\prime}(|\nabla u|+|\nabla v_r|) |\nabla u - \nabla v_r|^2\, \mathrm{d}x \leq c \left( \omega(2r)^\frac{p}{q} + r^{\min\{\beta,\gamma\}} \right) \left( \widetilde{\varphi}\left(\dashint_{B_{2r}(x_0)} \left|\nabla u \right| \, \mathrm{d}x \right) + \Lambda \right) \,,
%\label{eq:stimaV}
%\end{equation}
%for a constant $c=c(d,p,q,L,\kappa_0)$.
Here we follow the argument of \cite[Lemma 6.2]{HastoOk}. First, since $v_r$ is a weak solution to \eqref{eq:systemphitilde} and $u-v_r\in W^{1,\widetilde{\varphi}}_0(B_r(x_0))$ is an admissible test function for the weak formulation of \eqref{eq:systemphitilde}, we obtain
\begin{equation*}
\int_{B_r(x_0)} \frac{\widetilde{\varphi}^{\prime}(|\nabla v_r|)}{|\nabla v_r|}  \nabla v_r \cdot \nabla (u-v_r)\, \mathrm{d}x =0 \,.
\end{equation*} 
Now, with this and %\eqref{eq:equivtphi2} and 
\eqref{eq:equivtphi3}, we get
\begin{equation}
\begin{split}
%\int_{B_r(x_0)} \left| \mathbf{V}(\nabla u) - \mathbf{V}(\nabla v_r) \right|^2 \, dx & \sim   
\int_{B_r(x_0)} \widetilde{\varphi}^{\prime\prime}(|\nabla u|+|\nabla v_r|)& |\nabla u - \nabla v_r|^2\, \mathrm{d}x  \lesssim \int_{B_r(x_0)} \widetilde{\varphi}(|\nabla u|)\, \mathrm{d}x - \int_{B_r(x_0)} \widetilde{\varphi}(|\nabla v_r|)\, \mathrm{d}x \\
&- \int_{B_r(x_0)} \widetilde{\varphi}^{\prime}(|\nabla v_r|)  \frac{\nabla v_r}{|\nabla v_r|} \cdot (\nabla u - \nabla v_r)\, \mathrm{d}x \\
& = \int_{B_r(x_0)} \widetilde{\varphi}(|\nabla u|)\, \mathrm{d}x - \int_{B_r(x_0)} \widetilde{\varphi}(|\nabla v_r|)\, \mathrm{d}x \\
& = \int_{B_r(x_0)} \widetilde{\varphi}(|\nabla u|) - \varphi(x,|\nabla u|)\, \mathrm{d}x + \int_{B_r(x_0)} \varphi(x,|\nabla u|) - \varphi(x,|\nabla v_r|)\, \mathrm{d}x \\
& + \int_{B_r(x_0)} \varphi(x,|\nabla v_r|) - \widetilde{\varphi}(|\nabla v_r|) \, \mathrm{d}x =: I_1 + I_2 + I_3 \,.
\end{split}
\label{eq:stima1}
\end{equation}
We proceed to estimate each integral separately. 

Both the terms $I_1$ and $I_3$ can be estimated as 
\begin{equation}
\dashint_{B_r(x_0)} |\varphi(x,|\nabla u|)-\widetilde{\varphi}(|\nabla u|)|\,\mathrm{d}x\,,\,\,\dashint_{B_r(x_0)} |\varphi(x,|\nabla v_r|)-\widetilde{\varphi}(|\nabla v_r|)|\,\mathrm{d}x \leq c \left( \omega(2r)^\frac{p}{q} + r^\gamma \right) \left[\widetilde{\varphi}\left(\dashint_{B_{2r}(x_0)} |\nabla u|\,\mathrm{d}x\right) +\Lambda\right]
\label{eq:stimehok}
\end{equation}
by the very same argument of \cite[Lemma 6.2]{HastoOk} exploiting \azero{}, \inc{p}, \dec{q}, \eqref{eq:comptildephi}, \eqref{eq:5.16HOK} and \eqref{eq:5.17HOK}, where the constant $c$ depends on $d,p,q,L$. 

As for $I_2$, using $v_r$ in the definition of $u$ as almost minimizer of $\mathcal{F}$, and taking into account \eqref{eq:5.17HOK} and \eqref{eq:5.16HOK}, we get
\begin{equation*}
\begin{split}
I_2 & \leq \mathcal{L}^d(B_r) \kappa r^\beta \left( \dashint_{B_r(x_0)} \varphi(x,|\nabla v_r|) \, \mathrm{d}x + \lambda \right) +  \mathcal{L}^d(B_r) \lambda \\
& \leq c \kappa r^\beta \mathcal{L}^d(B_r) \widetilde{\varphi}\left (\dashint_{B_{2r}(x_0)}|\nabla {u}|\, \mathrm{d}x \right) + c \mathcal{L}^d(B_r) \Lambda \,. 
\end{split}
\end{equation*}
Plugging this estimate together with \eqref{eq:stimehok} into \eqref{eq:stima1} we obtain \eqref{eq:stimaV}. %\\
\end{proof}

\subsection{Local H\"{older} continuity}\label{sec:locholder}

In this section, we will  establish two main regularity results for almost minimizers, which will be instrumental in the proof of the local Lipschitz continuity result. First, we will prove the $C^{0,\alpha}$-regularity of any almost-minimizer $u$ for $\mathcal{F}$, locally within $\Omega$, for \emph{any} exponent $\alpha \in (0,1)$. Then, we will show that $\nabla u$ is locally $C^{0,\alpha}$ for \emph{some} exponent $\alpha \in (0,1)$ away from the free-boundary $\partial\{u>0\}$. 

We are now in position to prove the first regularity result. 
\begin{theorem}[$C^{0,\alpha}$-regularity]
Let $\varphi \in \Phic(\Omega)$, $\varphi(x,\cdot)\in C^1([0,\infty))$ be satisfying \vauno{}, and such that $\varphi_t$ comply with \azero{}, \inc{p-1} and \dec{q-1} for some $1<p\leq q$. If $u\in W^{1,\varphi}_{\rm loc}(\Omega)$ is an almost minimizer of $\mathcal{F}$ with constant $\kappa\leq \kappa_0$ and exponent $\beta$, then $u\in C^{0,\alpha}_{\rm loc}(\Omega)$ for any $\alpha\in(0,1)$. More precisely, if $\Omega^\prime\Subset\Omega$ is fixed, there exists $0<R_0<\frac{1}{2}{\rm dist}(\Omega^\prime, \partial \Omega)$ complying with \eqref{eq:smallradius} (thus depending on $u$) and a constant $c=c(d, p, q, \alpha, \kappa_0, R_0, \Omega')$ such that
\begin{equation}
[u]_{C^{\alpha}(\Omega^\prime)} \leq c \left(\int_\Omega|\nabla u|\,\mathrm{d}x + \Lambda 
\right)\,. 
\label{eq:c0alphabound}
\end{equation}
\label{thm:C0alpha}
\end{theorem}
\begin{proof}
Let $R_0\in(0,1)$ be sufficiently small to be determined later. Let $\Omega^\prime\Subset\Omega$ be fixed, and {assume that \eqref{eq:smallradius} holds for $R_0$, $R_0<{\rm dist}(\Omega^\prime, \partial \Omega)/2$.} For any $B_{2r}\Subset\Omega^\prime$, with $2r\leq R_0$, denote by $v_r$ the $\widetilde{\varphi}$-harmonic replacement of $u$ in $B_r$. Set ${\eta}(r):= \omega(2r)^\frac{p}{q} + r^{\min\{\beta,\gamma\}}$, and note that ${\eta}( \cdot )\leq 2$.

Let $\tau\in(0,\frac{1}{2})$. Then, by \eqref{eq:equivtphi5} with $\varepsilon=\frac{1}{2}$, $z_1=\nabla u$ and $z_2=\nabla v_r$, Proposition \ref{lem:stimaV}, \eqref{eq:5.15DSV} for $w=v_r$, \eqref{eq:comptildephi} %\eqref{eq:equivtphi4}
 and \eqref{eq:5.17HOK} %the definition of $v_r$ 
we get
\begin{equation}
\begin{split}
\int_{B_{\tau r}} \widetilde{\varphi}(|\nabla u|)\,\mathrm{d}x & \lesssim  \int_{B_r}\widetilde{\varphi}^{\prime\prime}(|\nabla u|+|\nabla v_r|) |\nabla u - \nabla v_r|^2\, \mathrm{d}x + \int_{B_{\tau r}} \widetilde{\varphi}(|\nabla v_r|)\,\mathrm{d}x \\
& \leq c  r^d {\eta}(2r)  \widetilde{\varphi}\left(\dashint_{B_{2r}(x_0)} \left|\nabla u \right| \, \mathrm{d}x\right) + c \Lambda r^d + c(\tau r)^d \sup_{B_{r/2}} \widetilde{\varphi}(|\nabla v_r|) \\
& \leq c  r^d {\eta}(R_0)  \widetilde{\varphi}\left(\dashint_{B_{2r}(x_0)} \left|\nabla u \right| \, \mathrm{d}x\right) + c \Lambda r^d + c(\tau r)^d \dashint_{B_r} \widetilde{\varphi}( |\nabla v_r|) \,\mathrm{d}x \\
& \leq c  r^d {\eta}(R_0)  \widetilde{\varphi}\left(\dashint_{B_{2r}} \left|\nabla u \right| \, \mathrm{d}x\right) + c \Lambda r^d + c(\tau r)^d \dashint_{B_r} ({\varphi}\left( x,|\nabla v_r|\right)+1 ) \,\mathrm{d}x \\
& \leq c  r^d {\eta}(R_0)  \widetilde{\varphi}\left(\dashint_{B_{2r}} \left|\nabla u \right| \, \mathrm{d}x\right) + c \Lambda r^d + c(\tau r)^d \widetilde{\varphi}\left(\dashint_{B_{2r}} \left|\nabla u \right| \, \mathrm{d}x\right) + c \Lambda(\tau r)^d  \\
& \leq c  r^d \left({\eta}(R_0)  + \tau^d\right) \widetilde{\varphi}\left(\dashint_{B_{2r}} \left|\nabla u \right| \, \mathrm{d}x\right) + c \Lambda r^d \,,
\end{split}
\end{equation}
whence, by Jensen's inequality, 
\begin{equation}
\begin{split}
\widetilde{\varphi}\left(\dashint_{B_{\tau r}}|\nabla u|\,\mathrm{d}x \right) \leq \dashint_{B_{\tau r}} \widetilde{\varphi}(|\nabla u|)\,\mathrm{d}x &  \leq c  \left(\frac{{\eta}(R_0)}{\tau^d}  + 1 \right) \widetilde{\varphi}\left(\dashint_{B_{2r}} \left|\nabla u \right| \, \mathrm{d}x\right) + c \frac{\Lambda}{\tau^d} \\
& \leq c \widetilde{\varphi}\left( \left(\frac{{\eta}(R_0)}{\tau^d}  + 1 \right) \dashint_{B_{2r}} \left|\nabla u \right| \, \mathrm{d}x + \frac{\Lambda}{\tau^d} \right)\,,
\end{split}
\end{equation}
so that
\begin{equation}
\int_{B_{\tau r}}|\nabla u|\,\mathrm{d}x \lesssim (\eta(R_0)+\tau^d) \int_{B_{2r}} |\nabla u|\,\mathrm{d}x + \Lambda r^d\,.
\end{equation}

The previous estimate trivially holds also for $\tau\in(\frac{1}{2},2)$. Let $\sigma\in(0,d)$  and choose $R_0$ small enough such that
\begin{equation}
{\eta}(R_0) \leq \varepsilon_0\,,
\end{equation}
where $\varepsilon_0$ is that of Lemma \ref{lem:iterationlemma},  applied to the function $f(\tau):=\int_{B_{\tau r}} |\nabla u|\,\mathrm{d}x$. We then obtain the Morrey-type estimate
\begin{equation}
\dashint_{B_{\rho}} |\nabla u|\,\mathrm{d}x \leq c   \left(\left(\frac{\rho}{R_0}\right)^{-\sigma} \dashint_{B_{R_0}} |\nabla u|\,\mathrm{d}x + \Lambda \rho^{-\sigma}\right) \quad \mbox{ for all balls $B_{\rho}\subset \Omega^\prime$ with $\rho\in(0,R_0]$,}
\label{eq:intMorrey}
\end{equation}
where the constant $c\geq1$ depends on $d,p,q,L,\kappa_0$ and $\sigma$.

Now, let $\alpha\in(0,1)$ be arbitrarily fixed. Then choosing $\sigma=1-\alpha$ in \eqref{eq:intMorrey}, by Morrey's Theorem (see, e.g., \cite[Chapter III, Theorem 1.1]{GIAQUINTA}) we infer that $u\in C^{0,\alpha}(\overline{B_{\rho}})$, for every $B_\rho\subset B_{R_0}$ and
\begin{equation}
[u]_{C^{\alpha}(B_{\rho})} \leq c \left[ R_0^{1-\alpha}\, \dashint_{B_{R_0}} |\nabla u|\,\mathrm{d}x + \Lambda \right]\,. 
\label{eq:holderonballs}
\end{equation}

Therefore, $u\in C^{0,\alpha}_{\rm loc}(\Omega)$ for every $\alpha\in(0,1)$ fixed. Since $\overline{\Omega'}$ is compact, \eqref{eq:c0alphabound} follows from \eqref{eq:holderonballs} by a standard covering argument. 
\end{proof}

With the following result, we establish the local H\"{o}lder continuity of the gradient of almost-minimizers, away from the free boundary.
\begin{theorem}[$C^{1,\alpha}$-regularity in $\{u>0\}$]
Let $\varphi \in \Phic(\Omega)$, $\varphi(x,\cdot)\in C^1([0,\infty))$ be such that $\varphi_t$ comply with \azero{}, \inc{p-1} and \dec{q-1} for some $1<p\leq q$. If $u\in W^{1,\varphi}_{\rm loc}(\Omega)$ is an almost minimizer of $\mathcal{F}$ with constant $\kappa\leq \kappa_0$ and exponent $\beta$ and $\varphi$ satisfies \vauno{} with 
\begin{equation}
\omega(r)\leq c r^\theta, \quad \mbox{ for all $r\in(0,1]$ and for some $\theta\in(0,1)$,}
\label{eq:omegar}
\end{equation}
then $u\in C^{1,\bar{\alpha}}_{\rm loc}(\{u>0\})$ for some $\bar{\alpha}\in(0,1)$ depending on $d,p,q,L,\beta,\theta$. More precisely, for any $\widetilde{\Omega}\Subset \{u>0\}\cap\Omega$, there exists an exponent $\bar{\alpha}=\bar{\alpha}(d,p,q,L,\beta,\theta)$ and a constant $C=C(\widetilde{\Omega},d,p,q,L,\kappa_0,\beta,\theta)$ such that 
\begin{equation}
[\nabla u]_{C^{\bar{\alpha}}(\widetilde{\Omega})} \leq C \left(\int_{\Omega} |\nabla u|\,\mathrm{d}x + \Lambda \right)\,.
\label{eq:c1alphabound}
\end{equation} 
\label{thm:c1alphabound}
\end{theorem}
\begin{proof}
Let $\widetilde{\Omega}\Subset \{u>0\}\cap\Omega$. Then, from \eqref{eq:intMorrey} and a standard covering argument we deduce that, for every $\sigma\in(0,d)$,
\begin{equation}
\dashint_{B_{2r}} |\nabla u|\,\mathrm{d}x \leq   c_\sigma   \left( \int_{\Omega} |\nabla u|\,\mathrm{d}x + \Lambda \right) r^{-\sigma} \quad \mbox{ for all balls $B_{2r}\subset \widetilde{\Omega}$ with $2r\in(0,R_0]$,}
\label{eq:(12)}
\end{equation}
where $R_0$ is that of Theorem \ref{thm:C0alpha} and $c_\sigma=c_\sigma(d,p,q,L,R_0,\sigma)$. 

Let us fix any of such balls, say $B_r\Subset \widetilde{\Omega}$, with $r\leq R_0/2$ small enough and to be determined later, and denote by $v_r$ the $\widetilde{\varphi}$-harmonic replacement of $u$ in $B_r$.  Arguing as in the estimate \eqref{eq:stima1} we get
\begin{equation}
\begin{split}
\dashint_{B_{r}} \widetilde{\varphi}^{\prime\prime}(|\nabla u|+|\nabla v_r|) |\nabla u - \nabla v_r|^2\, \mathrm{d}x & \lesssim \dashint_{B_{r}} \widetilde{\varphi}(|\nabla u|) - \varphi(x,|\nabla u|)\, \mathrm{d}x + \dashint_{B_{r}} \varphi(x,|\nabla u|) - \varphi(x,|\nabla v_r|)\, \mathrm{d}x \\
& \,\,\,\,\,\, + \dashint_{B_{ r}} \varphi(x,|\nabla v_r|) - \widetilde{\varphi}(|\nabla v_r|) \, \mathrm{d}x \,.
\end{split}
\label{eq:estim1}
\end{equation}
Since, by the convex-hull property, $v_r(B_r)$ is contained in the convex hull of $v_r(\partial B_r)=u(\partial B_r)$ and so in particular $\chi_{\{u>0\}}= \chi_{\{v_r>0\}}$ on $B_r$, taking into account also \eqref{eq:stimehok} , \eqref{eq:omegar}, \eqref{eq:comptildephi} and \eqref{eq:5.16HOK} we get
\begin{equation}
\begin{split}
\dashint_{B_r} \varphi(x,|\nabla u|) - \varphi(x,|\nabla v_r|)\, \mathrm{d}x & \leq \kappa r^\beta  \dashint_{B_r} \varphi(x,|\nabla v_r|)\, \mathrm{d}x + \kappa r^\beta\lambda \dashint_{B_r}\chi_{\{v_r>0\}}(x)\,\mathrm{d}x \\
& \leq \kappa r^\beta  \dashint_{B_r} (\varphi(x,|\nabla v_r|)-\widetilde{\varphi}(|\nabla v_r|))\, \mathrm{d}x + \kappa r^\beta \dashint_{B_r} \widetilde{\varphi}(|\nabla {u}|))\, \mathrm{d}x \\
& \,\,\,\,\,\, + \lambda\kappa r^{\beta}  \\
& \leq \kappa r^\beta \left(r^\frac{\theta p}{q} + r^\gamma \right) \left[\widetilde{\varphi}\left(\dashint_{B_{2r}}|\nabla u| \,\mathrm{d}x \right)+\Lambda \right] + \kappa r^\beta \dashint_{B_r} {\varphi}(x, |\nabla {u}|))\, \mathrm{d}x \\ 
& \,\,\,\,\,\, + \Lambda\kappa r^{\beta}  \\
& = \kappa r^\beta \left(r^\frac{\theta p}{q} + r^\gamma + 1 \right) \left[\widetilde{\varphi}\left(\dashint_{B_{2r}} |\nabla u|\,\mathrm{d}x \right) + \Lambda \right] \\
& \leq c \kappa r^\beta \left[\widetilde{\varphi}\left(\dashint_{B_{2r}} |\nabla u|\,\mathrm{d}x \right) + \Lambda \right]  \,.
\end{split}
\label{eq:estim2}
\end{equation}
Now, plugging \eqref{eq:estim2} into \eqref{eq:estim1}, and estimating the other two integrals again as in \eqref{eq:stimehok}, we obtain
\begin{equation}
\begin{split}
\dashint_{B_r} \widetilde{\varphi}^{\prime\prime}(|\nabla u|+|\nabla v_r|) |\nabla u - \nabla v_r|^2\, \mathrm{d}x & \lesssim (\kappa_0 r^\beta + r^{\gamma_1}   ) \left[\widetilde{\varphi}\left(\dashint_{B_{2r}} |\nabla u|\,\mathrm{d}x \right) + \Lambda \right]  \,,
\end{split}
\label{eq:estim3}
\end{equation}
where $\gamma_1:=\min\{\frac{\theta p}{q}, \gamma\}<1$. %From \eqref{eq:estim3}, arguing as for the proof of \eqref{eq:stimasenzaV} we get
Following the argument in \cite[Lemma 6.3]{HastoOk}, %However, we prefer to write down the details since the following computations will be useful also in the rest of the paper. 
we set $\eta(r):= r^\beta + r^{\gamma_1}\leq 2$, and apply \eqref{eq:equivtphi5} with 
$\varepsilon=\sqrt{\eta(r)}$, \eqref{eq:comptildephi}, Lemma \ref{lem:lemma5.15HOK}, \eqref{eq:estim3} and \eqref{eq:stimaV}. We find that
\begin{align*}
\dashint_{B_r}&\widetilde{\varphi}(|\nabla u-\nabla v_r|)\,\mathrm{d}x \\
&\lesssim \sqrt{\eta(r)} \dashint_{B_r}\left[\widetilde{\varphi}(|\nabla u|)+\widetilde{\varphi}(|\nabla v_r|)\right] \,\mathrm{d}x + \frac{1}{\sqrt{\eta(r)}}\dashint_{B_r}\widetilde{\varphi}''(|\nabla u|+|\nabla v_r|)|\nabla u- \nabla v_r|^2\,\mathrm{d}x \\
&\lesssim 
\sqrt{\eta(r)} \dashint_{B_r}\left[\varphi(x,|\nabla u|) + \varphi(x,|\nabla v_r|)+1\right]\,\mathrm{d}x
+ \sqrt{\eta(r)}\left(\widetilde{\varphi}\left(\fint_{B_{2r}}|\nabla u|\,\mathrm{d}x \right)+\Lambda\right)\\
&\lesssim \sqrt{\eta(r)}\left(\widetilde{\varphi}\left(\dashint_{B_{2r}}|\nabla u|\,\mathrm{d}x \right)+\Lambda\right).
\end{align*}
Therefore, by Jensen's inequality and \dec{q} of $\widetilde{\varphi}$, we have 
\[
\widetilde{\varphi}\left(\dashint_{B_r}|\nabla u-\nabla v_r|\, \mathrm{d}x\right)
\leq \dashint_{B_r}\widetilde{\varphi}(|\nabla u-\nabla v_r|)\, \mathrm{d}x
\lesssim \widetilde{\varphi}\left(\eta(r)^{\frac{1}{2q}}\left(\dashint_{B_{2r}}|\nabla u|\,\mathrm{d}x+\Lambda\right)\right),
\]
whence, since $\widetilde{\varphi}$ is strictly increasing, we finally obtain
\begin{equation}
\begin{split}
\dashint_{B_r} \left|\nabla u - \nabla v_r \right| \, dx & \lesssim r^{\beta_1} \left(\dashint_{B_{2r}} |\nabla u|\,\mathrm{d}x  + \Lambda \right)  \,,
\end{split}
\label{eq:estim3bis}
\end{equation}
where $\beta_1:=\min\{\frac{\beta}{2q},\frac{\gamma_1}{2q}\}<1$ and the implicit constant depends also on $\kappa_0$. 

On the other hand, with $\tau\in (0,\frac{1}{2})$, 
\begin{equation}
\begin{split}
\dashint_{B_{\tau r}} \left|\nabla u - (\nabla u)_{B_{\tau r}} \right| \, \mathrm{d}x & \leq c \dashint_{B_{\tau r}} \left| \nabla u - (\nabla v_r)_{B_{\tau r}} \right| \, \mathrm{d}x \\
& \leq c \tau^{-d} \dashint_{B_{r}} \left|\nabla u - \nabla v_r \right|\, \mathrm{d}x + c \dashint_{B_{\tau r}} \left|\nabla v_r - (\nabla v_r)_{B_{\tau r}} \right|\, \mathrm{d}x \,.
\end{split}
\label{eq:estim4}
\end{equation}
Note that
%\begin{equation}
%\begin{split}
%\left| (\mathbf{V}(\nabla {u}^{*}_r))_{B_{\tau r}(x_0)} - (\mathbf{V}(\nabla {u}))_{B_{\tau r}(x_0)} \right|^2 & = \left| \dashint_{B_{\tau r}(x_0)}(\mathbf{V}(\nabla {u}^{*}_r) - \mathbf{V}(\nabla {u})) \,\mathrm{d}x  \right|^2 \\
%& \leq \dashint_{B_{\tau r}(x_0)}\left|\mathbf{V}(\nabla {u}^{*}_r) - \mathbf{V}(\nabla {u})\right |^2\,\mathrm{d}x \\
%& \leq \tau^{-d} \dashint_{B_{r}(x_0)} \left| \mathbf{V}(\nabla {u}) - \mathbf{V}(\nabla {u^*_r}) \right|^2 \, \mathrm{d}x \,,
%\end{split}
%\label{eq:estim5}
%\end{equation}
using \eqref{eq:DSVestim} and \eqref{eq:5.18HOK}  we obtain %and $u^*_r$, 
%and the minimality of $(\nabla v_r)_{B_{r}}$
\begin{equation}
\begin{split}
\dashint_{B_{\tau r}} \left|\nabla v_r - (\nabla v_r)_{B_{\tau r}} \right|\, \mathrm{d}x  & \leq c \tau^{\mu_0} \dashint_{B_{r/2}} \left|\nabla v_r\right| \, \mathrm{d}x \\
%c \tau^{\mu_0} \dashint_{B_{r/2}} \left| \nabla v_r - (\nabla v_r)_{B_{r/2}} \right| \, \mathrm{d}x \\
%& \leq c \tau^{\mu_0} \dashint_{B_{r/2}} \left|\nabla v_r\right| \, \mathrm{d}x \\
& \leq c \tau^{\mu_0} \left(\dashint_{B_{r}} |\nabla u|\,\mathrm{d}x + \Lambda\right)\,. 
\end{split}
\label{eq:estim6}
\end{equation}
Inserting the estimate \eqref{eq:estim6} into \eqref{eq:estim4}, and taking into account \eqref{eq:estim3bis}, we get
\begin{equation}
\begin{split}
\dashint_{B_{\tau r}} \left|\nabla u - (\nabla u)_{B_{\tau r}} \right| \, \mathrm{d}x & \leq c \tau^{-d} \dashint_{B_{r}} \left|\nabla u - \nabla v_r \right|\, \mathrm{d}x + c \tau^{\mu_0} \left(\dashint_{B_{2r}} |\nabla u|\,\mathrm{d}x + \Lambda\right) \\
%& \leq c \left [\tau^{-d} (\kappa r^\beta + r^{\gamma_1}) +\tau^\mu \right ]  \left(\dashint_{B_{2r}} |\nabla u|\,\mathrm{d}x + \Lambda\right) \\
& \leq c \left [\tau^{-d} r^{\beta_1 } +\tau^{\mu_0} \right ] \left(\dashint_{B_{2r}} |\nabla u|\,\mathrm{d}x + \Lambda\right) \,.
\end{split}
\label{eq:estim7}
\end{equation}
%where $\beta_1:=\min\{\beta,\gamma_1\}<1$. 

Now, with \eqref{eq:(12)}, we get
\begin{equation}
\begin{split}
\dashint_{B_{\tau r}} \left|\nabla u - (\nabla u)_{B_{\tau r}} \right| \, \mathrm{d}x \leq c \left [\tau^{-d} r^{\beta_1 -\sigma} +\tau^{\mu_0} r^{-\sigma} \right ] \left(\int_{\Omega} |\nabla u|\,\mathrm{d}x + \Lambda \right) \,. 
\end{split}
\label{eq:estim8}
\end{equation}
Choosing $\tau:= r^{\frac{\beta_1}{{\mu_0} + d}}$, we have $\tau^{-d} r^{\beta_1 -\sigma} = \tau^{\mu_0} r^{-\sigma}= r^{\frac{\beta_1 {\mu_0}}{{\mu_0}+d}-\sigma}$, so that
\begin{equation}
\begin{split}
\dashint_{B_{\tau r}} \left|\nabla u - (\nabla u)_{B_{\tau r}} \right| \, \mathrm{d}x & \leq c  r^{\frac{\beta_1 \mu_0}{\mu_0+d}-\sigma} \left(\int_{\Omega} |\nabla u|\,\mathrm{d}x + \Lambda \right) \\
& = c  (\tau r)^{\frac{\beta_1 {\mu_0} - \sigma({\mu_0} + d)}{{\mu_0} + d + \beta_1}} \left(\int_{\Omega} |\nabla u|\,\mathrm{d}x + \Lambda \right) \,, 
\end{split}
\label{eq:estim9}
\end{equation}
and if $r\leq\frac{1}{2}\min\{R_0, 4^{-\frac{{\mu_0} + d}{\beta_1}}\}$, then $\rho:=\tau r < \frac{r}{2}$. Fixed $\bar{\sigma}:= \frac{{\mu_0} \beta_1}{2({\mu_0}+d)}$ and setting $\bar{\alpha}:= \frac{\beta_1 {\mu_0} - \bar{\sigma}(\mu + d)}{{\mu_0} + d + \beta_1}=\frac{\beta_1 {\mu_0}}{2({\mu_0} + d + \beta_1)}<1$, 
 the previous estimate \eqref{eq:estim9} can be rewritten as
\begin{equation}
\dashint_{B_{\rho}} \left|\nabla u - (\nabla u)_{B_{\rho}} \right| \, \mathrm{d}x  \leq  c   \left(\int_{\Omega} |\nabla u|\,\mathrm{d}x + \Lambda \right) \rho^{\bar{\alpha}} \,, 
\label{eq:estim10}
\end{equation}
where $B_{\rho}\subset B_{2r}\subset \widetilde{\Omega}$. From the Campanato-type embedding (see, e.g., \cite[Chapter III, Theorem 1.3]{GIAQUINTA}), this implies $\nabla {u}\in C^{0,{{\bar{\alpha}}}}_{\rm loc}(\widetilde{\Omega})$ and, since $\widetilde{\Omega}$ was arbitrary, $\nabla {u}\in C^{0,{{\bar{\alpha}}}}_{\rm loc}(\{u>0\})$. %Now, recalling that $\mathbf{V}^{-1}$ is $C^{0,\gamma}(\mathbb{R}^d)$ for some $\gamma=\gamma(p,q)$ (see \cite[Lemma 2.10]{DSV}), we conclude that $\nabla u\in C^{0,\gamma{{\bar{\alpha}}}}_{\rm loc}(\{u>0\})$. 
Moreover, \eqref{eq:c1alphabound} follows by a covering argument as in the end of the proof of  Theorem \ref{thm:c1alphabound}.  
%\begin{equation*}
%[\nabla u]_{C^{\bar{\alpha}}(\widetilde{\Omega})} \leq C \left(\int_{\Omega} |\nabla u|\,\mathrm{d}x + \Lambda \right)\,.
%\end{equation*}
\end{proof}

%\begin{lemma}
%Let $(\varphi_j)_{j\in\mathbb{N}}$ be a sequence of $\Phic(\Omega)\cap C^1([0,+\infty))$, such that $(\varphi_j)_t$ satisfies \azero{}, \inc{p-1}, \dec{q-1} uniformly, where $1<p\le q<+\infty$, and $\varphi_j$ satisfies \vauno{} uniformly. Let $(v_j)$ be a sequence of [...] functions. Assume that
%\begin{enumerate}
%\item[(i)] \begin{equation}\varphi_j(y,t)\to \varphi_\infty(t) \,\, \mbox{ uniformly on $B_1\times K$, where {$K\subset [0,+\infty)$} is compact, }%|\xi|^{\bar p}+h_\infty(x_0,\xi)=:f_\infty(\xi).
%\label{eq:unifconvphij}
%\end{equation}
%for a $C^1([0,+\infty))$ convex function $\varphi_\infty$ such that $\varphi'_\infty$ satisfies \inc{p-1} and \dec{q-1};
%\item[(ii)] $(v_j)$ is equibounded in $L^\infty(\Omega)$ and complying with
%\begin{equation}
%\int_{B_R} \varphi_j(x,|\nabla v_j|)\,\mathrm{d}x \leq \int_{B_R} \varphi_j(x,|\nabla w|)\,\mathrm{d}x + o(1)\,, \quad \mbox{ as $j\to+\infty$,}
%\end{equation}
%for every $w\in W^{1,\infty}(B_R)$;
%\item[(iii)] there exists $v_0\in W^{1,1}(B_R)$ such that $v_j\rightharpoonup v_0$ weakly in $W^{1,p}(B_R)$ and such that
%\begin{equation}
%\int_{B_R} \varphi_{\infty}(|\nabla v_0|)\,\mathrm{d}x \leq \mathop{\lim\inf}_{j\to+\infty} \int_{B_R} \varphi_j(x,|\nabla v_j|)\,\mathrm{d}x \,.
%\end{equation}
%\end{enumerate}
%Then $v_0$ is $\varphi_\infty$-harmonic in $B_R$. 
%\end{lemma}

\subsection{Proof of the local Lipschitz continuity} \label{sec:lipschitzcontinuous}

\begin{lemma}
Under the assumptions of Theorem \ref{thm:c1alphabound}, let $u$ be a bounded almost minimizer of $\mathcal{F}$ in $B_1(0)$ with constant $\kappa\leq \kappa_0$ and exponent $\beta$. Assume that $B_1(0)=\{u>0\}$. Then
\begin{equation}
|\nabla u(0)| \leq  C\,,
\label{eq:boundgrad0}
\end{equation}  
where the constant $C$ depends on  $p,q,d,L,\kappa_0,\beta, \theta, \Lambda, \|u\|_{L^\infty(B_1(0))}$. 
\label{lem:boundgrad0}
\end{lemma}

\proof
%Clearly, it will suffice to prove \eqref{eq:boundgrad0} in the case that $\|u\|_{L^\infty(B_1(0))}<+\infty$, as otherwise there is nothing to prove.

We have
\begin{equation}
|\nabla u(0)| \leq |\nabla u(0)-(\nabla u)_{B_{\frac{1}{4}}(0)}| + |(\nabla u)_{B_{\frac{1}{4}}(0)}|\,.
\label{eq:triangular1}
\end{equation}
Now, the first term on the right hand side above can be estimated by Theorem \ref{thm:c1alphabound}: observe that, since we are assuming $u$ to be bounded, by Remark \ref{rmk:smallness}, (ii) the constant appearing there depends on  $p,q,d,L,\kappa_0,\beta,  \theta, \Lambda, \|u\|_{L^\infty(B_1(0))}$}. We then have 
\begin{equation}
\begin{split}
|\nabla u(0)-(\nabla u)_{B_{\frac{1}{4}}(0)}| & \leq \dashint_{B_{\frac{1}{4}}(0)}|\nabla u(0)-\nabla u(x)|\,\mathrm{d}x \\
& \leq 4^{-\alpha}[\nabla u]_{C^{\alpha}(B_{\frac{1}{4}}(0))} \\
&    \leq C \left(\int_{B_{\frac{1}{2}}(0)}|\nabla u| \,\mathrm{d}x + \Lambda \right)\,,
\end{split}
\label{eq:triangular2}
\end{equation}
while
\begin{equation}
\begin{split}
|(\nabla u)_{B_{\frac{1}{4}}(0)}| & \leq  \mathcal{L}^d(B_\frac{1}{4}) \left( \int_{B_{\frac{1}{2}}(0)} |\nabla u|\,\mathrm{d}x + \Lambda \right) \,.
\end{split} 
\label{eq:triangular3}
\end{equation}
Further, from the Caccioppoli inequality Lemma \ref{lem:caccioppoli}, %\eqref{eq:comptildephi}, 
\eqref{v0phi}, \inc{1} and the boundedness of $u$ on $B_1(0)$, we get
\begin{equation}
\begin{split}
\int_{B_{\frac{1}{2}}(0)} |\nabla u|\,\mathrm{d}x \leq Lq\int_{B_{\frac{1}{2}}(0)} \left(\varphi(x,|\nabla u|)+1\right)\,\mathrm{d}x  & \leq c  \int_{B_{1}(0)} \left(\varphi(x,2\|u\|_{L^\infty(B_1(0))})+\Lambda\right)\,\mathrm{d}x\\
& \leq c \left(\max \left\{\|u\|^p_{L^\infty(B_1(0))}, \|u\|^q_{L^\infty(B_1(0))}\right\} + \Lambda \right)\,.
\end{split}
\label{eq:triangular4}
\end{equation}

Combining \eqref{eq:triangular1}--\eqref{eq:triangular4} we obtain \eqref{eq:boundgrad0}, and this concludes the proof. 

\endproof

\begin{lemma}
Let $R>0$ be such that $B_{2R}(0)\Subset \Omega$, and $(r_j)_{j\in\mathbb{N}}$, $(\sigma_j)_{j\in \mathbb{N}}$  be sequences of nonnegative numbers, with $R<\frac{1}{2r_j}$ for every $j$, $r_j\to0$ as $j\to+\infty$, and 
\begin{equation} 
\sigma_j\to+\infty \quad \mbox{ and } \quad \varphi(0,\sigma_j) r_j\to0 \,\,\,\, \mbox{ as $j\to+\infty$}. 
\label{eq:sigmaj}
\end{equation}
We define, for every $j$,
\begin{equation}
\varphi_j(x,t):= \frac{\varphi(r_j x, \sigma_j t)}{\varphi(0,\sigma_j)}\,, \,\, x\in B_{2R}(0)\,,\,\, t>0\,.
\label{eq:phij}
\end{equation} 
Then, 
\begin{enumerate}
\item[(i)] the functions
\begin{equation}
\varphi^-_j(t):= \inf_{y\in B_{2R}(0)} \varphi_j(y,t)\,, \quad \mbox{ and } \quad \varphi^+_j(t):= \sup_{y\in B_{2R}(0)} \varphi_j(y,t)\,,
\end{equation}
are weak $\Phi$ functions satisfying \inc{p} and \dec{q}. Moreover, for $j$ large enough, 
\begin{equation}
\min\{t^p,t^q\}\lesssim \varphi_j^-(t)\leq \max\{t^p,t^q\}, \quad \min\{t^p,t^q\}\leq \varphi_j^+(t)\lesssim \max\{t^p,t^q\},
\label{eq:stimevarphij}
\end{equation}
where the hidden constants are independent of $j$; 
\item[(ii)]  there exists $j_0\in\mathbb{N}$ such that $\varphi_j$ complies with \azero\ for $j\geq j_0$ with $L=2$;
\item[(iii)] there exists $j_0\in\mathbb{N}$ such that $\varphi_j$ complies with \vauno{} for $j\geq j_0$ with the same $\omega$;
\item[(iv)] there exists a convex function $\varphi_\infty \in C^1([0,+\infty))$, whose derivative $\varphi'_\infty$ complies with \inc{p-1} and \dec{q-1} such that
\begin{equation}
\varphi_j(x,t)\to \varphi_\infty(t) \,\, \mbox{ uniformly on $B_{2R}(0)\times K$, where {$K\subset [0,+\infty)$} is compact. }
\label{eq:unifconvphij}
\end{equation} 
\end{enumerate}
\label{lem:asymptphi}
\end{lemma}
\begin{proof} 
For the proof of $(i)$ we can argue as in \cite[pp. 13--14]{LSSV1}. Concerning the proof of $(ii)$, 
let $j_0\in\mathbb{N}$ be such that
\begin{equation}
r_j<1\,,\,\, \sigma_j>L^{\frac 1p}\,,\,\, \varphi(0,\sigma_j)>1\,, \,\, \mbox{ and } \,\, \varphi(0,\sigma_j) r_j <\min\{1,1/{\rm diam}(\Omega)\gamma_d\} \quad \mbox{ for $j\geq j_0$}.
\label{eq:smallnesss}
\end{equation}
Thanks to \eqref{eq:smallnesss} it is easy to show that $\varphi_{2r_jR}^-(\sigma_j)\in[\omega(2r_jR),\frac{1}{\mathcal{L}^d(B_{2r_jR})}]$. Then applying \vauno \  for $\varphi$ and using  that $\omega(2r_jR)\le 1$, we deduce $(ii)$.

We turn to the proof of $(iii)$.
Let $j_0\in\mathbb{N}$ be such that \eqref{eq:smallnesss} holds.
We show that $\varphi_j$ satisfies \vauno{} in $B_{2R}(0)$ for $j\geq j_0$. Let $j$ as above be fixed and let $\tau\in(0,1)$. We have to prove that
\begin{equation}
\sup_{x\in B_{2\tau R}(0)}\varphi_j(x,t)\le (1+\omega(2\tau R))\inf_{x\in B_{2\tau R}(0)}\varphi_j(x,t), \quad\forall t>0 \,\, \mbox{ s. t. } \,\,\inf_{x\in B_{2\tau R}(0)}\varphi_j(x,t)\in \left [\omega(2\tau R), \frac{1}{\mathcal{L}^d(B_{2\tau R})} \right].
\label{eq:claimva1}
\end{equation} 
Note that
\begin{equation}
\inf_{x\in B_{2\tau R}(0)}\varphi_j(x,t)\in \left [\omega(2\tau R), \frac{1}{\mathcal{L}^d(B_{2\tau R})} \right] \iff \varphi^-_{B_{2\tau R r_j}(0)}(\sigma_j t) \in \left [\varphi(0,\sigma_j)\omega(2\tau R), \frac{\varphi(0,\sigma_j)}{\mathcal{L}^d(B_{2\tau R})} \right]\,.
\end{equation}
Now, from \eqref{eq:smallnesss}, 
\begin{equation}
\left [\varphi(0,\sigma_j)\omega(2\tau R), \frac{\varphi(0,\sigma_j)}{\mathcal{L}^d(B_{2\tau R})} \right] \subseteq \left [\omega(2\tau R r_j), \frac{\varphi(0,\sigma_j)r_j^d}{\mathcal{L}^d(B_{2\tau R r_j})} \right] \subseteq \left [\omega(2\tau R r_j), \frac{1}{\mathcal{L}^d(B_{2\tau R r_j})} \right] \,.
\end{equation}
Then, by \vauno{} for $\varphi$, $r_j<1$ and the fact that $\omega$ is increasing, we get
\begin{equation}
 \varphi^+_{B_{2\tau R r_j}(0)}(\sigma_j t)\le (1+\omega(2\tau R)) \varphi^-_{B_{2\tau R r_j}(0)}(\sigma_j t), \quad\forall t>0 \,\, \mbox{ s. t. } \,\,\inf_{x\in B_{2\tau R}(0)}\varphi_j(x,t)\in \left [\omega(2\tau R), \frac{1}{\mathcal{L}^d(B_{2\tau R})} \right]\,,
\end{equation} 
whence \eqref{eq:claimva1} follows up to dividing both the sides by $\varphi(0,\sigma_j)$.  
The proof of assertion \eqref{eq:unifconvphij} is postponed to the Appendix, see Lemma \ref{lem:asympphij}. 
\end{proof}

Let $(\varphi_j)_{j\in\mathbb{N}}$ be the sequence defined in \eqref{eq:phij} and, correspondingly, consider the scaled functional
\begin{equation}
\hat{\mathcal{F}}_j(v,\Omega):= \int_\Omega \varphi_j(x,|\nabla v|)\,\mathrm{d}x + \frac{\lambda}{\varphi(0, \sigma_j)}  \int_\Omega \chi_{\{v>0\}}(x)\,\mathrm{d}x \,.
\label{eq:scaledfunct}
\end{equation}
With given $u$, we also consider for every $j$ the blow-up function
\begin{equation}
v_j(x):= \frac{u(r_j x)}{\sigma_j r_j}\,, \quad x\in B_{2R}(0)\,.
\label{eq:blowupfunct}
\end{equation}
We then have the following result about the asymptotic behavior of a blow-up sequence defined by scaling an almost minimizer of $\mathcal{F}$.

\begin{proposition}\label{lem:scaledv}
Let $R$, $(r_j)_{j\in\mathbb{N}}$, $(\sigma_j)_{j\in\mathbb{N}}$, $(\varphi_j)_{j\in\mathbb{N}}$ and $\varphi_\infty$ be as in Lemma \ref{lem:asymptphi}. %Let $0<R<\frac{1}{r_j}$, with $r_j\to0$ as $j\to+\infty$, and $\sigma_j>0$ be such that
%\begin{equation} 
%\sigma_j\to+\infty \quad \mbox{ and } \quad \varphi(0,\sigma_j) r_j\to0 \,\,\,\, \mbox{ as $j\to+\infty$}. 
%\label{eq:sigmaj}
%\end{equation}
Let $u$ be a bounded almost minimizer of $\mathcal{F}$ in $B_2(0)$ with constant $\kappa\leq \kappa_0$ and exponent $\beta$. 
%Set $\psi(t):=\varphi(0,t)$ and define %$\gamma_j>0$ such that $\gamma_j\to+\infty$ and $\gamma_j r_j\to0$ as $j\to+\infty$. Set $\psi(t):=\varphi(0,t)$ and define 
%\begin{equation}
%\sigma_j:=\psi^{-1}({\textstyle\frac{1}{\gamma_j r_j}}). 
%\label{eq:sigmaj}
%\end{equation}
Then, for every $j$, the function $v_j$ defined in \eqref{eq:blowupfunct}
%\begin{equation}
%v_j(x):= \frac{u(r_j x)}{\sigma_j r_j}\,, \quad x\in B_{2R}(0)\,,
%\end{equation}
is an almost minimizer of the scaled functional $\hat{\mathcal{F}}_j$ \eqref{eq:scaledfunct} in $B_{2R}(0)$,
%\begin{equation}
%\hat{\mathcal{F}}_j(v,\Omega):= \int_\Omega \varphi_j(x,|\nabla v|)\,\mathrm{d}x + \frac{\lambda}{\varphi(0, \sigma_j)}  \int_\Omega \chi_{\{v>0\}}(x)\,\mathrm{d}x \,,
%\label{eq:scaledfunct}
%\end{equation}
with constant $\hat{\kappa}:=\kappa r_j^\beta$ and the same exponent $\beta$. 
%   where
%\textcolor{blue}{\begin{equation}
%\varphi_j(x,t):= \frac{\varphi(r_j x, \sigma_j t)}{\varphi(0,\sigma_j)}\,, \,\, x\in B_{2R}(0)\,,\,\, t>0\,. 
%\label{eq:phij}
%\end{equation} }
Moreover, 
%\begin{enumerate}
%\item[(i)] there exists a convex function $\varphi_\infty \in C^1([0,+\infty))$, whose derivative $\varphi'_\infty$ complies with \inc{p-1} and \dec{q-1} such that
%\begin{equation}
%\varphi_j(x,t)\to \varphi_\infty(t) \,\, \mbox{ uniformly on $B_{2R}(0)\times K$, where {$K\subset [0,+\infty)$} is compact; }
%\label{eq:unifconvphij}
%\end{equation} 
%\item[(ii)] 
if $\|v_j\|_{L^\infty(B_{2R}(0))}\leq M$, there exists $v_\infty\in W^{1,1}(B_R(0))$ such that, up to a subsequence, $v_j\rightharpoonup v_\infty$ weakly in $W^{1,p}(B_R(0))$, and uniformly in $B_R(0)$, and $v_\infty$ is $\varphi_\infty$-harmonic in $B_R(0)$.  
%\end{enumerate}
\label{lem:lemmascaled}
\end{proposition}

\begin{proof}
Let $B_\rho(x_0)$ be a ball such that $\overline{B_\rho(x_0)}\subset B_{\frac{1}{r_j}}(0)$, and $w\in W^{1,\varphi}(B_\rho(x_0))$ such that $w=v_j$ on $\partial B_\rho(x_0)$. Setting $y_0:=r_j x_0$, we then have
\begin{equation*}
u(y) = \sigma_j r_j w({\textstyle\frac{y}{r_j}})=:\widetilde{w}_j(y)\,, \quad \mbox{ on $\partial B_{r_j\rho}(y_0)$}
\end{equation*}
and, by the almost minimality of $u$, we get 
\begin{equation}
\int_{B_{r_j\rho}(y_0)} \varphi(y,|\nabla u(y)|)+ \lambda \chi_{\{{u} >0\}} (y) \, \mathrm{d}y \leq (1+\kappa(r_j\rho)^\beta)\int_{B_{r_j\rho}(y_0)} \varphi(y,|\nabla \widetilde{w}_j(y) |)+ \lambda \chi_{\{{\widetilde{w}_j} >0\}} (y) \, \mathrm{d}y\,.
\label{eq:almu}
\end{equation}
Now, with the change of variables $x=\frac{y}{r_j}$, we have
\begin{equation}
\begin{split}
\int_{B_{r_j\rho}(y_0)} \varphi(y,|\nabla \widetilde{w}_j(y) |)+ \lambda \chi_{\{{\widetilde{w}_j} >0\}} (y) \, \mathrm{d}y & = r_j^d \int_{B_{\rho}(x_0)} \varphi(r_j x,|\nabla \widetilde{w}_j(r_j x) |)+ \lambda \chi_{\{{\widetilde{w}_j} >0\}} (r_j x) \, \mathrm{d}x \\
& = r_j^d \int_{B_{\rho}(x_0)} \varphi(r_j x,\sigma_j |\nabla{w}(x)| )+ \lambda \chi_{\{{{w}} >0\}} (x) \, \mathrm{d}x
\end{split}
\label{eq:scaledw}
\end{equation}
and, in a similar way, 
\begin{equation}
\int_{B_{r_j\rho}(y_0)} \varphi(y,|\nabla u(y)| )+ \lambda \chi_{\{{u} >0\}} (y) \, \mathrm{d}y = r_j^d \int_{B_{\rho}(x_0)} \varphi(r_j x,\sigma_j |\nabla{v}_j(x))|+ \lambda \chi_{\{{v_j} >0\}} (x) \, \mathrm{d}x \,.
\label{eq:scaledu}
\end{equation}
Plugging \eqref{eq:scaledu} and \eqref{eq:scaledw} into \eqref{eq:almu}, and multiplying both the sides of the inequality by $\frac{1}{\varphi(0,\sigma_j)}$, and recalling the definition of $\varphi_j$, we then infer that $v_j$ is an almost minimizer of the functional $\hat{\mathcal{F}}_j$ defined in \eqref{eq:scaledfunct}. 

%{As for assertion $(i)$, it is an immediate consequence of Lemma \ref{lem:asympphij}, so we turn to the proof of $(ii)$.} 

Now, we notice that {applying} Lemma \ref{lem:caccioppoli} {to $\varphi_j$} we obtain for $v_j$ the Caccioppoli-type estimate 
\begin{equation}
\int_{B_{\rho}(y)} \varphi_j(x,|\nabla v_j|)\,\mathrm{d}x  \leq c \left( \int_{B_{2\rho}(y)} \varphi_j\left(x,\frac{|v_j-(v_j)_{y,2\rho}|}{2\rho}\right)\,\mathrm{d}x + { \frac{\lambda}{\varphi(0,\sigma_j)} \rho^d} \right) 
\label{eq:caccioppolitype}
\end{equation}
for any $B_{2\rho}(y)\Subset B_{2R}(0)$, where the constant $c$ only depends on $d,p,q,\kappa,\beta$, and is a uniform constant with respect to $j$. % Setting
%\begin{equation}
%\varphi^-_j(t):= \inf_{y\in B_{2R}(0)} \varphi_j(y,t)\,, \quad \mbox{ and } \quad \varphi^+_j(t):= \sup_{y\in B_{2R}(0)} \varphi_j(y,t)\,,
%\end{equation}
%and arguing as in \cite[pp. 13--14]{LSSV1}, it can be shown that both $\varphi_j^-$ and $\varphi_j^+$ are weak $\Phi$ functions satisfying \inc{p} and \dec{q}. Moreover, for $j$ large enough, 
%\begin{equation}
%\min\{t^p,t^q\}\lesssim \varphi_j^-(t)\leq \max\{t^p,t^q\}, \quad \min\{t^p,t^q\}\leq \varphi_j^+(t)\lesssim \max\{t^p,t^q\},
%\label{eq:stimevarphij}
%\end{equation}
%where the hidden constants are independent of $j$. 

%By assumption \vauno{},
%\begin{equation}\label{a01}
%\varphi_{j}^-(1)\ge\frac{1}{2} \varphi_{j}^+(1)\ge \frac{1}{2}
%\end{equation}
%if $\varphi_{r_j}^-(\sigma_j)\in[\omega(r_j), \frac{1}{\mathcal{L}^d(B_{r_j})}]$. For $j$ large enough,  $\varphi_{r_j}^-(\sigma_j)<\omega(r_j)\le 1$ does not occur since,  
%by \eqref{cons2} and \eqref{v0phi}, this would entail $\sigma_j$ equibounded. If in the end  $\varphi_{r_j}^-(\sigma_j)>\frac{1}{\mathcal{L}^d(B_{r_j})}$, then
%\begin{equation}\label{a02}
%\varphi_{j}^-(1) > 1,
%\end{equation}
%for $j$ large enough.

Recall that $\|v_j\|_{L^\infty(B_{2R}(0))}\leq M$. By \eqref{eq:stimevarphij}, $\varphi_j$ satisfy \inc{p} and \dec{q}, with constants independent of $j$. They also satisfy \azero{} and \vauno{} with $L$ and $\omega$ independent of $j$, by Lemma \ref{lem:asymptphi}. Then, the radius $r_0$ of Remark~\ref{rmk:smallness}(ii) can be chosen independently of $j$. It follows that, applying Proposition \ref{thm:C0alpha}  to each $v_j$, each of them is locally $\alpha$-H\"older continuous  on $B_{2R}(0)$, and the $C^{0,\alpha}$-estimate \eqref{eq:c0alphabound} on $B_R(0)$, holds with a uniform bound not depending on $j$.

Since $\varphi^-_j$ is \inc{p}, combining the bounds $\|v_j\|_{L^\infty(B_{2R}(0))}\leq M$, \eqref{eq:caccioppolitype} and \eqref{eq:stimevarphij} we obtain

%Then, arguing as in Remark~\ref{rmk:smallness}(ii) and using that $\|v_j\|_{L^\infty(B_{2R}(0))}\leq M$, \eqref{eq:stimevarphij} and $0<\frac{1}{\varphi(0,\sigma_j)}<1$ for $j$ large enough, %from \eqref{eq:caccioppolitype} 
%we then obtain an analogous estimate as \eqref{eq:newargument}, with $\Omega''=B_R(0)$ and $\Omega'=B_{2R}(0)$,
%\begin{equation}
%\begin{split}
%\int_{B_{R}(0)} \varphi_j(x,|\nabla v_j|)\,\mathrm{d}x  & \leq R^{-ds_0}\left(\int_{B_R(0)} %%\varphi_j(x,|\nabla v_j|)^{1+s_0}\,\mathrm{d}x\right)^\frac{1}{1+s_0}\leq \tilde{c} = \tilde{c}(d,p,q,L,\kappa_0,M,R,\lambda)%c \left( \int_{B_{2R}(0)} \varphi_j^+\left({\textstyle\frac{2M}{R}}\right)\,\mathrm{d}x + {\lambda R^d} \right) \\
% & \leq c R^d \left(\max\{({\textstyle\frac{2M}{R}})^p, ({\textstyle\frac{2M}{R}})^q\} + \lambda \right)
%\end{split}
%\label{eq:bound1}
%\end{equation}
%for $j$ large enough. Observe that $\varphi_j$ satisfy \azero{} and \vauno{} with $L$ and $\omega$ independent of $j$, by Lemma \ref{lem:asymptphi}. As $\|v_j\|_{L^\infty(B_{2R}(0))}\leq M$, then, the radius $r_0$ of Remark~\ref{rmk:smallness}(ii) can be chosen independently of $j$. It follows that, applying Proposition \ref{thm:C0alpha}  to each $v_j$, each of them is locally $\alpha$-H\"older continuous  on $B_{2R}(0)$, and the $C^{0,\alpha}$-estimate \eqref{eq:c0alphabound} on $B_R(0)$, holds with a uniform bound not depending on $j$.
%
%Since $\varphi^-_j$ is \inc{p}, \eqref{eq:bound1} and \eqref{eq:stimevarphij} imply that
\begin{equation}
\sup_{j\geq j_0}\int_{B_{\rho}(0)} |\nabla v_j|^p\,\mathrm{d}x \leq \sup_{j\geq j_0}\int_{B_{\rho}(0)} \left( \frac{\varphi^-_j(|\nabla v_j|)}{\varphi_j^-(1)}+1 \right) \,\mathrm{d}x \leq C\,,
\end{equation}
for a constant $C$ depending on $d,p,q,L,\kappa_0,M,R,\lambda$. %$d,p,q,\kappa,\beta,R,M$,
Hence the above inequality also holds for $\rho=R$, whence we infer the existence of a function $v_\infty\in W^{1,p}(B_R(0))$ such that, up to a subsequence,
\begin{equation}
v_j \rightharpoonup v_\infty \quad \mbox{ weakly in $W^{1,p}(B_R(0))$. }
\label{eq:weakconverg}
\end{equation}

By \eqref{eq:c0alphabound} on $B_R(0)$,  \ref{eq:weakconverg} also gives %(up to change the subsequence) to
\begin{equation*}
v_j \to v_\infty \quad \mbox{ uniformly in $B_R(0)$}
\end{equation*}
%for any $\alpha\in(0,1)$, %by using the $C^{0,\alpha}$-estimate \eqref{eq:c0alphabound} 
since the sequence $(v_j)$ is equibounded by $M$ on $B_R(0)$. 

So we are left to prove that $v_\infty$ is $\varphi_\infty$-harmonic in $B_R (0)$. We first notice that using \eqref{eq:caccioppolitype} for $\frac{R}{2}<\rho<R$, exploiting the uniform bound $\|v_j\|_{L^\infty(B_{2R}(0))}\leq M$ and letting $\rho\to R$, we also get that the sequence of positive measures $\mu_j:= \varphi_j(\cdot,|\nabla v_j|)\,\mathcal{L}^d$ is equibounded on $B_R(0)$. Thus, we can find a Radon measure $\mu$ on $B_R(0)$ such that $$\mu_j \rightharpoonup^* \mu \mbox{\,\,\,\, on $B_R(0)$}$$
up to a subsequence (not relabeled). Now, 
let us fix  $w\in W^{1,\varphi_\infty}(B_R(0))$ be such that $\{w\neq v_\infty\}\Subset B_R(0)$. Since $\varphi_\infty$ satisfies \dec{q}, we can find a sequence $(w^\varepsilon)_{\varepsilon>0}\subset W^{1,\infty}(B_R(0))$ of regularizations of $w$, strongly converging to $w$  in $W^{1,{\varphi_\infty}}(B_R(0))$ as $\varepsilon\to0$ (see, e.g., \cite[Lemma 6.4.5]{HH}). 

Let $\rho<\rho'\in (0, R)$, with %$\rho'$ such that \eqref{anche} holds, 
$\mu(\partial B_{\rho'})=\mu(\partial B_\rho)=0$ and
$\{w\neq v_\infty\}\Subset B_{\rho}$. Let $\eta\in C^\infty_c(B_{\rho'})$ be such that $\eta = 1$ on $B_\rho$, { $0\leq \eta\leq1$, $|\nabla \eta|\leq \frac{2}{\rho'-\rho}$,} and define
$\zeta_j= \eta w^\varepsilon+ (1-\eta) {v}_j$. Since $\{\zeta_j\neq {v}_j\}\Subset B_{\rho'}$, using the almost minimality of $v_j$, straightforward computations lead to
\begin{equation}\label{compar}
\begin{split}
&\int_{B_{\rho'}}  \varphi_j(x,|\nabla v_j|)\,\mathrm{d}x \leq (1+\kappa(r_j\rho')^\beta)\hat{\mathcal{F}}_j(\zeta_j,B_{\rho'})\\
& \leq \int_{B_{\rho}} \varphi_j(x,|\nabla w^\varepsilon|)\,\mathrm{d}x + c(1+\kappa(r_j\rho')^\beta) \int_{B_{\rho'}\setminus B_\rho}\left(\varphi_j(x,|\nabla v_j|)+\varphi_j(x,|\nabla w^\varepsilon|)+ %\varphi_j(y,|\nabla {\hat v}_j|)+
\varphi_j\left(x,\frac{|w^\varepsilon-{v}_j|}{\rho'-\rho}\right)\right)\mathrm{d}x\\
& \,\,\,\,\,\, +\frac{\lambda}{\varphi(0,\sigma_j)}\mathcal{L}^d(B_{\rho}) + \kappa(r_j\rho')^\beta \hat{\mathcal{F}}_j(w^\varepsilon,B_{\rho}) + (1+\kappa(r_j\rho')^\beta) \frac{\lambda}{\varphi(0,\sigma_j)}\, \mathcal{L}^d(B_{\rho'}\setminus B_\rho)
\end{split}
\end{equation}
for a suitable constant $c\ge 1$ depending only on $L$ and $p, q$. First, we note that  
\begin{equation*}
\displaystyle\mathop{\lim\sup}_{j\to+\infty}\int_{B_{\rho'}\setminus B_\rho}\varphi_j(x,|\nabla v_j|)\,\mathrm{d}x \leq \mu(B_{\rho'}\setminus B_\rho)\quad \mbox{ and } \quad \hat{\mathcal{F}}_j(w^\varepsilon,B_{\rho}) \leq C \mathcal{L}^d(B_\rho)
\end{equation*}
for $j$ sufficiently large, %First, we note that  H\"{o}lder inequality and \eqref{eq:bound1} give
%\begin{equation*}
%\int_{B_{\rho'}\setminus B_\rho}\varphi_j(x,|\nabla v_j|)\,\mathrm{d}x \leq c\,\mathcal{L}^d(B_{\rho'}\setminus B_\rho)^{\frac{s_0}{1+s_0}} \quad \mbox{ and } \quad \hat{\mathcal{F}}_j(w^\varepsilon,B_{\rho}) \leq C \mathcal{L}^d(B_\rho)
%\end{equation*}
%for $j$ sufficiently large, 
and
\begin{equation}
\frac{\lambda}{\varphi(0,\sigma_j)}\mathcal{L}^d(B_{\rho'}) + \kappa(r_j\rho')^\beta \hat{\mathcal{F}}_j(w^\varepsilon,B_{\rho}) + (1+\kappa(r_j\rho')^\beta) \frac{\lambda}{\varphi(0,\sigma_j)}\, \mathcal{L}^d(B_{\rho'}\setminus B_\rho) \to 0
\end{equation}
as $j\to+\infty$, for fixed $\rho,\rho', \varepsilon$.

Now we deal with the convergence of the integral terms above. Using the uniform convergence \eqref{eq:unifconvphij} %of $\varphi_j$ to $\varphi_\infty$ on $B_1\times K$, where $K$ is any compact subset in $\mathbb{R}$, 
we have that 
\[
\lim_{j\to+\infty}\int_{B_{\rho'}\setminus B_\rho}\varphi_j(x,|\nabla w^\varepsilon|)\,\mathrm{d} x=\int_{B_{\rho'}\setminus B_\rho}\varphi_\infty(|\nabla w^\varepsilon|)\,\mathrm{d}x,
\]
since $|\nabla w^\varepsilon|$ is bounded. Likewise, we have
%As for the other term, we first notice that $\varphi_j(\cdot,t)\leq\varphi_j^+(t)$,
 %the boundedness of $w^\varepsilon$ and $v_j$, and a similar argument as for \eqref{eq:bound1} %and \eqref{convergetronc}(ii) 
%entail the equi-integrability of $\left\{\varphi_j\left(\cdot,\frac{|w^\varepsilon-{v}_j|}{\rho'-\rho}\right)\right\}_{j\in\mathbb{N}}$. Furthermore, taking into account the  pointwise convergence of  $\varphi_j\left(x,\frac{|w^\varepsilon-{v}_j|}{\rho'-\rho}\right)$ to $\varphi_\infty\left(\frac{|w^\varepsilon-{v}_\infty|}{\rho'-\rho}\right)$ {implied by \eqref{eq:unifconvphij}} 
%we may appeal to Vitali convergence theorem, which ensures that
\[
\lim_{j\to+\infty}\int_{B_{\rho'}\setminus B_\rho}\varphi_j\left(x,\frac{|w^\varepsilon-{v}_j|}{\rho'-\rho}\right)\,\mathrm{d}x=\int_{B_{\rho'}\setminus B_\rho}\varphi_\infty\left(\frac{|w^\varepsilon-{v}_\infty|}{\rho'-\rho}\right)\,\mathrm{d}x\,.
\]
Therefore, passing to the liminf as $j\to+\infty$ in \eqref{compar}, we have
\[
\begin{split}
\mathop{\lim\inf}_{j\to+\infty} & \int_{B_{\rho'}} \varphi_j(x,|\nabla v_j|)\,\mathrm{d}x \\
& \leq \int_{B_\rho}\varphi_\infty(|\nabla w^\varepsilon|)\,\mathrm{d}x+c\left[\int_{B_{\rho'}\setminus B_\rho}\left(\varphi_\infty(|\nabla w^\varepsilon|)+\varphi_\infty\left(\frac{|w^\varepsilon-v_\infty|}{\rho'-\rho}\right)\right)\mathrm{d}x\right]+\mu(B_{\rho'}\setminus B_\rho)\,. %c\,\mathcal{L}^d(B_{\rho'}\setminus B_\rho)^{\frac{s_0}{1+s_0}}\,.
\end{split}
\]
Now we  let $\varepsilon\to0$ and, recalling that $w=v_\infty$ outside $B_\rho$, we easily obtain 
\[
\mathop{\lim\inf}_{j\to+\infty}\int_{B_{\rho}} \varphi_j(x,|\nabla v_j|)\,\mathrm{d}x\le \int_{B_\rho}\varphi_\infty(|\nabla w|)\,\mathrm{d}x+c\int_{B_{\rho'}\setminus B_\rho}\varphi_\infty(|\nabla w|)\mathrm{d}x+ \mu(B_{\rho'}\setminus B_\rho)\,.  %c\,\mathcal{L}^d(B_{\rho'}\setminus B_\rho)^{\frac{s_0}{1+s_0}}\,.
\]
Therefore, with the lower semicontinuity result \eqref{lsc}, letting $\rho'$ tend to $\rho$ we finally get that for %${\mathcal L}^1$ -a.e. 
every $\rho\in(0,R)$ and any $w\in W^{1,{\varphi_\infty}} (B_R)$ such that $\{w\neq v_\infty\}\Subset B_\rho$ we have 
\[
\int_{B_\rho}\varphi_\infty(|\nabla v_\infty|)\,\mathrm{d}x\leq \int_{B_\rho}\varphi_\infty(|\nabla w|)\,\mathrm{d}x\,,
\]
as desired. 
\end{proof}

In order to prove the Lipschitz continuity of an almost minimizer, a tool will be the following Proposition, where we show that a bounded almost minimizer of $\mathcal{F}$ is sublinear in a neighborhood of a free-boundary point. 

\begin{proposition}\label{prop:stimaM}
Let $u$ be an almost minimizer of $\mathcal{F}$ in $B_1(x_0)$, where $x_0 \in \partial\{u>0\}\cap\Omega$, such that
\begin{equation}
\sup_{x\in B_1(x_0)} u(x) \leq M\,.
\end{equation}
Then there exists a constant $C_0=C_0(d,p,q,L,\kappa_0,\beta)\geq1$ such that
\begin{equation}
0\leq u(x) \leq C_0M|x-x_0|
\label{eq:lingrowth}
\end{equation}
for all $x\in B_r(x_0)$ and any $0<r<1$. 
\label{prop:lip1}
\end{proposition}

\proof
We may assume, without loss of generality, that $x_0=0$, and, throughout the proof, we will omit the center in the notation for a ball centered at $x_0$. We set
\begin{equation}
S(k,u):= \sup_{x\in B_{r_k}}|u(x)|\,,\quad r_k:=2^{-k}\,,\quad k\geq0\,,
\end{equation}
and our aim is to prove that there exists a constant $C\geq1$ such that
\begin{equation}
S(k+1,u) \leq \max\left\{CM r_{k+1}, \frac{S(k,u)}{2}\right\}\,,\quad \mbox{ for every} \,\, k\geq0\,.
\label{eq:iteration1}
\end{equation}
%where $\psi(t):=\varphi(0,t)$.

Indeed, once \eqref{eq:iteration1} has been established, arguing by induction we can prove that
\begin{equation}
S(k,u) \leq C M r_k\,,\quad \mbox{ for every} \,\, k\geq0\,.
\label{eq:iteration2}
\end{equation}
From this, given $r\in(0,1]$ and chosen $k\geq0$ such that $r_{k+1} < r \leq r_k$, we obtain
\begin{equation}
\|u\|_{L^\infty(B_r)} \leq \|u\|_{L^\infty(B_{r_k})} = S(k,u) \leq C M r_k = 2CM r_{k+1} \leq 2C M r\,,
\end{equation}
and then \eqref{eq:lingrowth}, with $C_0:=2C\geq1$.

In order to prove \eqref{eq:iteration1}, we argue by contradiction, and, for every $j\geq1$, we assume the existence of $u_j$ almost minimizer of $\mathcal{F}$ in $B_1$, with constant $\kappa$ and exponent $\beta$, and of an integer $k_j$ such that
\begin{equation}
S(k_j+1,u_j) > \max\left\{{jM r_{k_j+1}}, \frac{S(k_j,u_j)}{2} \right\}\,,
\label{eq:iteration3}
\end{equation}
%where $S_{k_j}:= \sup_{x\in B_{r_{k_j}}}|u_j(x)|$. 
Note that $\|u_j\|_{L^\infty(B_1)}\leq M$ implies $k_j\to+\infty$, since by \eqref{eq:iteration3} we infer $k_j>\log_2(j)-1$ for every $j$. Furthermore, with the uniform bound $\|u_j\|_{L^\infty(B_1)}\leq M$ and the same argument used in Proposition~\ref{lem:scaledv}, we can show that for any $\eta\in(0,1)$ the $u_j$ are uniformly locally $\eta$-H\"{o}lder continuous in $B_1$. Since $u_j(0)=0$, we obtain %Moreover,
%exploiting the H\"{older} continuity of $u_j$, the fact that $u_j(0)=0$, $|u_j|\leq M$ and \eqref{eq:c0alphabound} combined with the Caccioppoli inequality Lemma \ref{lem:caccioppoli}, for any fixed $\eta\in(0,1)$ we get
\begin{equation}
\sup_{x\in B_{r_{k_j+1}}}|u_j(x)| = \sup_{x\in B_{r_{k_j+1}}}|u_j(x) - u_j(0)| \leq  C_\eta r_{k_j+1}^\eta < C_\eta r_{k_j}^\eta\,,
\label{eq:supdecay}
\end{equation}
where $C_\eta$ is independent of $j$.

Now, we set
\begin{equation}
\sigma_j:= \frac{S(k_j+1,u_j)}{r_{k_j}}\,, 
\end{equation}
and we consider the scaled function
\begin{equation}
v_j(x) := \frac{u_j(r_{k_j} x)}{\sigma_j r_{k_j}} = \frac{u_j(r_{k_j} x)}{S(k_j +1 ,u_j)} \,, \quad x\in B_{\frac{1}{r_{k_j}}}\,.
\label{eq:vj}
\end{equation}
Note that, by \eqref{eq:iteration3}, 
\begin{equation}
\sigma_j \geq  j{\textstyle \frac{M}{2}} \to +\infty \quad \mbox{as $j\to+\infty$.}
\label{eq:sigmajj}
\end{equation}
Setting
\begin{equation}
\varphi_j(x,t):= \frac{\varphi(r_{k_j} x ,\sigma_j t)}{\varphi(0,\sigma_j)}\,, 
\end{equation}
with this choice of $\sigma_j$ and $\varphi_j$ we introduce the scaled functional $\hat{\mathcal{F}}_j$ defined as in \eqref{eq:scaledfunct}. Since, by \eqref{eq:sigmajj}, $\sigma_j>1$ for $j$ large enough, and $\varphi(0,t)$ is \dec{q}, we have, in view of \eqref{eq:supdecay} for $\eta=1-\frac{1}{2q}$
\begin{equation*}
\varphi(0,\sigma_j)r_{k_j} \leq \varphi(0,1) \sigma_j^q r_{k_j} = \varphi(0,1) \left(\frac{S(k_j+1,u_j)}{r_{k_j}^{1-\frac{1}{2q}}}\right)^q r_{k_j}^\frac{1}{2} \leq \varphi(0,1) C_{1-\frac{1}{2q}}^q r_{k_j}^\frac{1}{2}\,, \quad \mbox{ for $j$ large enough, }
\end{equation*}
whence \eqref{eq:sigmaj} follows. 

Taking into account \eqref{eq:iteration3}, for $x\in B_1$, we have
\begin{equation}
v_j(x) \leq \frac{S(k_j,u_j)}{S(k_j+1,u_j)} \leq 2 \frac{S(k_j,u_j)}{S(k_j,u_j)} = 2.
\label{eq:vjbounds}
\end{equation}
By Proposition~\ref{lem:scaledv} $v_j$ is an almost minimizer of $\hat{\mathcal{F}}_j$ in $B_1(0)$ with constant $\kappa r_{k_j}^\beta\leq \kappa_0$ and exponent $\beta$, there exist a $C^1([0,+\infty))$ convex function $\varphi_\infty$ whose derivative $\varphi'_\infty$ complies with \inc{p-1} and \dec{q-1}, a function $v_\infty\in W^{1,1}(B_1(0))$ such that, up to a subsequence, $v_j \to v_\infty$ uniformly in $B_1(0)$, and $v_\infty$ is $\varphi_\infty$-harmonic in $B_\frac{1}{2}(0)$.  Since, by \eqref{eq:vjbounds}, it holds that $0\leq v_\infty\leq 2$ in $B_\frac{1}{2}(0)$, and $v_\infty(0)=0$ being $v_j(0)\equiv 0$, from the strong minimum principle we must have $v_\infty\equiv 0$ in $B_{\frac{1}{2}}(0)$. However, from \eqref{eq:vj} we deduce that $\displaystyle\sup_{x\in B_{\frac{1}{2}}(0)}v_\infty(x)=1$ and this gives a contradiction. The proof is concluded. 

\endproof

We are now in position to prove the main result, Theorem~\ref{T1}. 

\proof[Proof of Theorem~\ref{T1}] Let $u$ be an almost minimizer of $\mathcal{F}$ in $\Omega$, with constant $\kappa\leq\kappa_0$ and exponent $\beta$. Let $\widetilde{\Omega}\Subset \Omega$, define $r_0$ as in \eqref{eq:smallradius} and set
\begin{equation*}
r_1:= \frac{1}{4} \min\left\{2r_0, {\rm dist}(\widetilde{\Omega},\partial\Omega)\right\} \quad \mbox{ and } \quad \Omega_{r_1}:= \left\{x\in\Omega:\,\, {\rm dist}(x,\partial\Omega)\geq r_1\right\}. 
\end{equation*} 
We recall that, by virtue of Theorem \ref{thm:C0alpha}, $u\in C^{0,\alpha}(\Omega_{r_1})$ for any fixed $\alpha\in(0,1)$, and set $M:=\|u\|_{L^\infty(\Omega_{r_1})}$.  Observe that this value $M$ depends on $r_0$, and hence on $u$ and $\widetilde{\Omega}$ via the integral  
\[
\int_{\widetilde{\Omega}} \varphi(x, |\nabla u|)^{1+s_0}\,\mathrm{d}x.
\]

Now, let $x_0\in \widetilde{\Omega}\cap\{u>0\}$ be arbitrarily fixed and, in order to estimate $|\nabla u(x_0)|$ we distinguish between two cases, according to $\tau:={\rm dist}(x_0, \partial\{u>0\}\cap \Omega)$. 

Let $\tau\leq r_1$ first, and choose $y_0\in \partial\{u>0\}\cap\Omega$ such that $|y_0-x_0|=\tau$. Since $B_{2\tau}(y_0)\subset \Omega_{r_1}$, we have $|u|\leq M$ in $B_{2\tau}(y_0)$. Then, by virtue of Proposition \ref{prop:lip1}, for every $x\in B_\tau(x_0)\subset B_{2\tau}(y_0)$ we have
\begin{equation}
u(x) \leq CM |x-y_0| \leq 2CM\tau.
\label{eq:lip2}
\end{equation} 

Now, let us consider the scaled function $u_\tau(x):=\frac{u(x_0+\tau x)}{\tau}$, $x\in B_1(0)$. Since $u$ is an almost minimizer of $\mathcal{F}$ in $B_\tau(x_0)$ with constant $\kappa$ and exponent $\beta$, a simple computation shows that $u_\tau$ is an almost minimizer in $B_1(0)$, with constant $\kappa \tau^\beta$ and exponent $\beta$, of the functional $\mathcal{F}_\tau$ defined as
\begin{equation*}
\mathcal{F}_\tau (w, \Omega):= \int_\Omega \varphi_\tau(x,|\nabla w|)\,\mathrm{d}x + {\lambda}  \int_\Omega \chi_{\{w>0\}}(x)\,\mathrm{d}x\,,
\end{equation*}
where $\varphi_\tau(x,t):=\varphi(x_0+\tau x,t)$. It is easy to check that $\varphi_\tau \in \Phic(B_1(0))$, $\varphi_\tau(x,\cdot)\in C^1([0,\infty))$ and that $(\varphi_\tau)_t$ complies with \azero{}, \inc{p-1} and \dec{q-1}. We only have to remark that also \vauno{} holds. For this, let $\rho\in(0,1)$, $B_\rho(y)\subset B_1(0)$ and recall that $\varphi$ satisfies \vauno{} on $B_{\tau\rho}(x_0+\tau y)$, so that
\begin{equation*}
\varphi^+_{B_{\tau\rho}(x_0+\tau y)}(t) \leq (1+\omega(\tau\rho)) \varphi^-_{B_{\tau\rho}(x_0+\tau y)}(t)\,, \quad \mbox{ if } \varphi^-_{B_{\tau\rho}(x_0 + \tau y)}(t)\in\left[\omega(\tau\rho), 1/\mathcal{L}^d(B_{\tau\rho})\right]\,.
\end{equation*}
Now, since $\varphi_\tau^\pm(t) = \varphi^\pm_{B_{\tau\rho}(x_0+\tau y)}(t)$, where $\varphi_\tau^\pm(t)$ are computed on $B_\rho(y)$, and $\tau\leq1$, from the previous estimate we infer
\begin{equation*}
\varphi^+_\tau(t) \leq (1+\rho^\beta) \varphi^-_\tau(t)\,, \quad \mbox{ if } \varphi^-_\tau(t)\in\left[\rho^\beta, 1/\mathcal{L}^d(B_{\rho})\right]\,.
\end{equation*}

Moreover, by \eqref{eq:lip2}, $|u_\tau|\leq 2CM$ in $B_1(0)$. Therefore, by Lemma \ref{lem:boundgrad0}, we deduce that
\begin{equation*}
|\nabla u(x_0)| = |\nabla u_\tau(0)| \leq \widetilde{C}\,,
\end{equation*}
where the constant $\widetilde{C}$ depends on  $p,q,d, L, \kappa_0,\beta, \Lambda, u, \widetilde{\Omega}$.

If, instead, $\tau\geq r_1$, we can perform an analogous argument as before with $u_{r_1}(x):=\frac{u(x_0+ r_1 x)}{r_1}$, $x\in B_1(0)$ in place of $u_\tau$, which satisfies $\|u_{r_1}\|_{L^\infty(B_1(0))}\leq \frac{M}{r_1}$, and $\mathcal{F}_{r_1}$ in place of $\mathcal{F}_\tau$. This concludes the proof. 

\endproof

\appendix
\section{}
\label{sec:appendix}
%\addcontentsline{toc}{section}{Appendices}

We recall a technical result (see \cite[Lemma 2.19]{LSSV1}) about the $\varphi$-recession function associated to a sequence of convex functions $\varphi_j$, capturing the behaviour at infinity of $\varphi_j$. 
\begin{lemma}\label{lem:lemma3.2fmt}
Let $(\varphi_j)_{j\in\mathbb{N}}$, $\varphi_j: [0,+\infty)\to[0,+\infty)$, be a sequence of $C^1$ convex functions satisfying $\varphi_j(0)=0$ and assume that $\varphi'_j$ satisfies \inc{p-1} and \dec{q-1}, where $1<p\le q<+\infty$.
Let $(\beta_j)\subset(0,\infty)$ be a sequence such that $\lim_j \beta_j=+\infty$. Then, setting
\begin{equation*}
\overline{\varphi}_j(t):=\frac{\varphi_{j}(t\beta_{j})}{\varphi_j(\beta_j)}\,,\quad t\in[0,+\infty)\,,\,\, j\in\mathbb{N}\,,
\end{equation*}
there exists a subsequence $(\beta_{j_k})$ such that $\overline{\varphi}_{j_k}$ converge to a  $C^1$ convex function $\varphi_\infty$ uniformly on compact subsets of $[0,+\infty)$. Moreover, $\varphi'_\infty$ satisfies \inc{p-1} and \dec{q-1}.
\end{lemma}

We are now in position to prove Lemma \ref{lem:asymptphi}$(iv)$. 

\begin{lemma}
Let $\varphi_j$ be the sequence defined as in \eqref{eq:phij}. Then \eqref{eq:unifconvphij} holds.
%\begin{equation}
%\varphi_j(y,t)\to \varphi_\infty(t) \,\, \mbox{ uniformly on $B_R\times K$, where {$K\subset [0,+\infty)$} is compact, }
%\label{eq:unifconvphij}
%\end{equation}
Moreover, $\varphi_\infty$ is a {$C^1([0,+\infty))$} function such that $\varphi'_\infty$ satisfies \inc{p-1} and \dec{q-1}.
\label{lem:asympphij}
\end{lemma}

\proof

We can apply Lemma~\ref{lem:lemma3.2fmt} to the constant sequence $\varphi_j(t)\equiv\varphi(0,t)$, with $\beta_j:=\sigma_j$.
Setting for brevity $\psi_j(t):=\varphi\left(0,t \sigma_j\right)$, and defining
\begin{equation*}
\overline{\varphi}_j(t):= \frac{\psi_j(t)}{\psi_j(1)} = \frac{\varphi(0, t\sigma_j)}{\varphi(0,\sigma_j)} \,,
\label{eq:scaledftilde}
\end{equation*}
we then obtain that, up to a subsequence, $\overline{\varphi}_j$ converges to a  $C^1$ convex function $\varphi_\infty$ uniformly on compact subsets of $[0,+\infty)$, with $\varphi'_\infty$ satisfying \inc{p-1} and \dec{q-1}. 
%we can apply Lemma~\ref{lem:lemma3.2fmt} , with {$\beta_j:=\sigma_j$}
%where ${\tilde f}(\xi)=|\xi|^{\bar p}+h(x_0,\xi)$,

Now, in order to prove \eqref{eq:unifconvphij}, it will suffice to show that, { with fixed $\tau>0$, there exists $j_1\geq1$ such that} %if $y\in B_1$, $t\in[0, R]$, 
\begin{equation}\label{step1}
|\varphi_j(y,t)-{\overline{\varphi}}_j(t)|\le \eta_{j,\tau}, \quad  {\mbox{for every $(y,t)\in B_R\times [0,\tau]$, for $j\geq j_1$, }}
\end{equation}
for some 
$\eta_{j,\tau}$ which is infinitesimal as $j\to+\infty$. 

The proof of the estimate \eqref{step1} follows by minor adaptations of the argument for the Step 3 of the proof of \cite[Lemma 3.1]{LSSV1}. For the reader's convenience, we provide the details.

First, we observe that
\[
\begin{split}
|\varphi_j(y,t)-\overline{\varphi}_j(t)|&= \frac{1}{\psi_j(1)}\left|\varphi(r_jy,t\sigma_j)-\psi_j(t)\right| \le \omega(r_j)\frac{\varphi^-_{r_j}(t\sigma_j)}{\psi_j(1)}, 
\end{split}
\]
if $\varphi_{r_j}^-(t\sigma_j)\in [\omega(r_j),{\textstyle\frac{1}{\mathcal{L}^d(B_{r_j})}}]$, thanks to \vauno{}. Recalling the definitions of $\varphi_{r_j}^-$ and $\psi_j$, together with \eqref{cons2}, the last term can be estimated as
\[
\omega(r_j)\frac{\varphi^-_{r_j}(t\sigma_j)}{\psi_j(1)}\le \omega(r_j)\overline{\varphi}_j(t) \le{\max\{\tau^p,\tau^q\}\,\omega(r_j)} \,. 
\]
On the other hand, if $\varphi^-_{r_j}(t\sigma_j)<\omega(r_j)$, from \eqref{cons2} and \eqref{v0phi}, we deduce 
\[
{\min\{(t\sigma_j)^p, (t\sigma_j)^q\}\le Lq \omega(r_j)\leq 1, \quad \mbox{ for $j$ large enough, }}
\]
entailing
\[
{t\sigma_j\le (Lq \omega(r_j))^{\frac{1}{q}}.}
\]
Then
\[
{|\varphi_j(y,t)-\overline{\varphi}_j(t)|\le \frac{2}{\psi_j(1)}(Lq \omega(r_j))^{\frac{p}{q}}\frac{L}{p}\lesssim \frac{1}{\psi_j(1)}\,.  }
\]
Finally, case $\varphi^-_{r_j}(t\sigma_j)>\frac{1}{\mathcal{L}^d(B_{r_j})}$ cannot occur for $j$ large enough since, {taking into account \eqref{cons2} and \eqref{eq:sigmaj},} it  would lead to  
\[
\frac{1}{\mathcal{L}^d(B_{r_j})}<\psi_j(t)\le\max\{\tau^p,\tau^q\}\psi_j(1)\,, %\le\max\{\tau^p,\tau^q\}\frac{\varepsilon_j}{r_j},
\]
that is 
\[
{\frac{1}{r_j^{d-1}}< \max\{\tau^p,\tau^q\}\gamma_d\,\varphi(0,\sigma_j) r_j, }
\]
{which clearly would give a contradiction for $j$ large.}

{Therefore, \eqref{step1} is proven { with $\eta_{j,\tau}:= \max\left\{\max\{\tau^p,\tau^q\}\,\omega(r_j), \frac{2}{\psi_j(1)}(Lq \omega(r_j))^{\frac{p}{q}}\frac{L}{p} \right\}$ for $j$ large enough}. }
\endproof

In order to prove the lower semicontinuity result of Lemma \ref{lem:semicontinuity}, we also need the following definitions and results about the maximal operator in Orlicz spaces (see \cite[Section 4.3]{HH}). 

\begin{definition}\label{maximal}
Given an open set $\Omega\subseteq\mathbb{R}^d$ and $f\in L^1_{\rm loc}(\Omega)$, the \emph{(centered) Hardy-Littlewood maximal operator} is $Mf:\Omega\to[0,\infty]$ defined as
\begin{equation}
    Mf(x) := \sup_{\rho>0}\frac{1}{\mathcal{L}^d(B_\rho(x))}\int_{B_\rho(x)\cap \Omega}|f(y)|\,\mathrm{d}y.
\label{eq:maximop}
\end{equation}
\end{definition}

The following result of boundedness for the local maximal operator can be found, e.g., in \cite[Corollary 1.9]{LSSV1}. 
\begin{lemma}\label{C:MaxOpBnd}
Let $\varphi\in\Phiw$ satisfy \inc{p} and \dec{q}, with $1<p\le q<+\infty$. Then there exists ${C}={C}(\varphi^{-1}(1),d,p,q)$ such that
\[
\int_\Omega \varphi(Mf)\,\mathrm{d} x\le {C} \int_\Omega \varphi(|f|)\,\mathrm{d} x
\]
for every  $f\in L^{\varphi}(\Omega)$ satisfying $\displaystyle{\int_\Omega\varphi(|f|)\,\mathrm{d} x\le1}$.
\end{lemma}

{We conclude the list of auxiliary results with the following Lusin-type approximation result in $W^{1,\varphi}$, which can be inferred from \cite[Theorem 3.3]{DIESTROVER10}. Indeed, for a fixed ball $B$, the argument therein can be applied to $\varphi^-_B$ without requiring $u\in W^{1,\varphi}_0(B)$, as the null extension outside $B$ is not needed taking into account the boundedness of the restricted maximal operator, Lemma \ref{C:MaxOpBnd}. }
\begin{theorem}\label{Lusin}
Let $\varphi\in\Phiw$ satisfy \inc{p} and \dec{q}, with $1<p\le q<+\infty$, and let $B\subset\mathbb{R}^d$ be a ball. For every $u\in W^{1,\varphi}(B)$ and every $\tau > 0$ there exists a function $u_\tau : B\to\mathbb{R}$, $u_\tau \in W^{1,\infty}(B)$ satisfying {\rm Lip}$(u_\tau) \le c\,\tau$ with $c=c(d)$, such that $u_\tau = u$ in
$\{M |\nabla u|\le \tau\}$ and
\[
{\mathcal L}^d(\{M |\nabla u| > \tau\})\le  \frac{1}{\varphi_{B}^-(\tau)}\int_{\{M|\nabla u|>\tau\}}\varphi_{B}^-(M |\nabla u|)\,\mathrm{d}x\,,
\]
{where $M$ is introduced in Definition \ref{maximal}.}
\end{theorem}

\begin{lemma}
Let $\varphi_j$ be the sequence defined as in \eqref{eq:phij}, and $\varphi_\infty$ be the function in \eqref{eq:unifconvphij}. Let $({v}_j)_{j\in\mathbb{N}}\subset W^{1,1}(B_{R}(0))$ be such that 
\begin{equation}
\sup_{j\in\mathbb{N}} \int_{B_R(0)} \varphi_j(y,|\nabla {v}_j|)\,\mathrm{d}y \leq C\,,
\label{assv3}
\end{equation}
and $v_j \to v_0$ $a.e.$ in $B_R(0)$. Then
\begin{equation}\label{lsc}
\displaystyle{\int_{B_R}\varphi_\infty(|\nabla v_0|)\,\mathrm{d} y\le \liminf_{j\to+\infty}\int_{B_R}\varphi_j(y,|\nabla {v}_j|)\,\mathrm{d} y}.
\end{equation}
\label{lem:semicontinuity}
\end{lemma}

\proof
Thanks to {the bound in \eqref{assv3}}, we apply Lemma \ref{C:MaxOpBnd} to $\varphi_j^-$, which is a weak $\Phi$-function satisfying \inc{p} and \dec{q}, obtaining
\begin{equation}\label{maxvj}
\int_{B_R(0)}\varphi_j^-(M|\nabla {v}_j|)\,\mathrm{d}y\le C,
\end{equation}
having taken into account that, thanks to \eqref{eq:stimevarphij}, $(\varphi_j^-)^{-1}(1)\simeq 1$, and the hidden constants do not depend on $j$. 
%By assumption \vauno{},
%\begin{equation}\label{a01}
%\varphi_{j}^-(1)\ge\frac{1}{2} \varphi_{j}^+(1)\ge \frac{1}{2}
%\end{equation}
%if $\varphi_{r_j}^-(\sigma_j)\in[\omega(r_j), \frac{1}{\mathcal{L}^d(B_{r_j})}]$. For $j$ large enough,  $\varphi_{r_j}^-(\sigma_j)<\omega(r_j)\le 1$ does not occur since,  
%by \eqref{cons2} and \eqref{v0phi}, this would entail $\sigma_j$ equibounded. If in the end  $\varphi_{r_j}^-(\sigma_j)>\frac{1}{\mathcal{L}^d(B_{r_j})}$, then
%\begin{equation}\label{a02}
%\varphi_{j}^-(1) > 1,
%\end{equation}
%for $j$ large enough. 
{By Chacon's Biting Lemma (see, e.g., \cite[Lemma~5.32]{Ambrosio-Fusco-Pallara:2000}) there exist a sequence of Borel subsets $A_h$ of $B_R(0)$ such that $\mathcal{L}^d(A_h)\to0$ as $h\to+\infty$, and a (not relabelled) subsequence such that $(\varphi_j^-(M|\nabla {v}_j|)\chi_{B_R(0)\setminus A_h})_j$ is equintegrable for every $h\geq1$. }

Let $\tau>1$. Then, applying Theorem~\ref{Lusin} to ${v}_j$, we find $v_j^\tau:B_R(0)\to\mathbb{R}$ such that
\begin{equation}
{\rm Lip}(v_j^\tau)\le c\,\tau\quad \hbox{ and }\quad v_j^\tau={v}_j \hbox{ in }B_R(0)\setminus E_j^\tau,
\label{eq:liptrunc}
\end{equation}
where $E_j^\tau := \{M|\nabla{\hat v}_j|>\tau\}$ and, by Chebychev's inequality, 
\begin{equation}\label{ejlambda}
{\mathcal L}^d(E_j^\tau\setminus A)\le  \frac{1}{\varphi_j^-(\tau)}\int_{E_j^\tau\setminus A}\varphi_j^-(M |\nabla { v}_j|)\,\mathrm{d} y,
\end{equation}
for any Borel set $A\subset B_R(0)$.
Moreover, from \eqref{ejlambda} with $A=\emptyset$, \eqref{maxvj}, {\inc{p} for $\varphi_j^-$, and the fact that by \eqref{eq:stimevarphij}, $\varphi_{j}^-(1)\gtrsim 1$ for $j$ large enough,} we deduce
\begin{equation}\label{misura}
{\mathcal{L}^d(E_j^\tau)\le %\int_{\{M|\nabla {\hat v}_j|>\lambda\}}\varphi_j^-(\lambda)\,\d y\le 
\frac{1}{\varphi_j^-(1)\tau^p}\int_{E_j^\tau}\varphi_j^-(M|\nabla {v}_j|)\,\mathrm{d} y\le \frac{C}{\tau^p}\,, }
\end{equation}
for $j$ large enough. 

We compute
\[
\begin{split}
\int_{B_R(0)}&\varphi_j(y,|\nabla { v}_j|)\,\mathrm{d}y\ge\int_{B_R(0)\setminus(A_h\cup E_j^\tau)}\varphi_j(y,|\nabla v_j^\tau|)\,\mathrm{d}y=\int_{B_R(0)\setminus A_h}\varphi_j(y,|\nabla v_j^\tau|)\,\mathrm{d}y\\
&-\int_{E_j^\tau\setminus A_h}\varphi_j(y,|\nabla v_j^\tau|)\,\mathrm{d}y=\int_{B_R(0)\setminus A_h}\left[\varphi_j(y,|\nabla v_j^\tau|)-\varphi_\infty(|\nabla v_j^\tau|)\right]\mathrm{d}y\\
&+\int_{B_R(0)\setminus A_h}\varphi_\infty(|\nabla v_j^\tau|)\,\mathrm{d}y-\int_{E_j^\tau\setminus A_h}\varphi_j(y,|\nabla v_j^\tau|)\,\mathrm{d}y.
\end{split}
\]
Since the convergence \eqref{eq:unifconvphij} implies
\[
\lim_{j\to+\infty}\int_{B_R(0)\setminus A_h}\left[\varphi_j(y,|\nabla v_j^\tau|)-\varphi_\infty(|\nabla v_j^\tau|)\right]\,\mathrm{d}y = 0\,,
\]
passing to the liminf in the previous inequality we obtain
\begin{equation}\label{liminf}
\liminf_{j\to+\infty}\int_{B_R(0)}\varphi_j(y,|\nabla {v}_j|)\,\mathrm{d}y\ge \liminf_{j\to+\infty}\int_{B_R(0)\setminus A_h}\varphi_\infty(|\nabla v_j^\tau|)\,\mathrm{d}y-
\limsup_{j\to+\infty}\int_{E_j^\tau\setminus A_h}\varphi_j(y,|\nabla v_j^\tau|)\,\mathrm{d}y.
\end{equation}
We are first dealing with the second term. We have
\[
\int_{E_j^\tau\setminus A_h}\varphi_j(y,|\nabla v_j^\tau|)\,\mathrm{d}y\le \int_{E_j^\tau\setminus A_h}\varphi^+_j(|\nabla v_j^\tau|)\,\mathrm{d}y.
\]
{In $E_j^\tau\setminus A_h$ we distinguish between the points of $B_R(0)$ where $\varphi_{r_j}^-(|\nabla{v}_j^\tau|\sigma_j)\in[\omega(r_j),1/\mathcal{L}^d(B_{r_j})]$, denoting the corresponding set by $S_{j,\tau}^1$, and the points where that condition does not hold. We then define
\begin{equation*}
S_{j,\tau}^2:= \left\{\varphi_{r_j}^-(|\nabla{v}_j^\tau|\sigma_j)<\omega(r_j)\right\}\cap B_R(0) \quad \mbox{ and }\quad S_{j,\tau}^3:=\left\{\varphi_{r_j}^-(|\nabla{v}_j^\tau|\sigma_j)>1/\mathcal{L}^d(B_{r_j})\right\}\cap B_R(0).
\end{equation*}
The set $S_{j,\tau}^3$ has to be empty for $j$ sufficiently large, as otherwise, using \eqref{cons2} for any fixed point therein, the resulting inequality {$\varphi(0,\sigma_j){\tau^q}>\frac{1}{\gamma_dr_j^d}$ would imply $\frac{1}{\varphi(0,\sigma_j) r_j}$ %$\frac{\gamma_j}{r_j^{d-1}}$ 
uniformly bounded with respect to $j$. }}

\noindent In {$S_{j,\tau}^2$}, thanks to \eqref{cons2} and \eqref{v0phi}, $\min\{(|\nabla{v}_j^\tau|\sigma_j)^p,(|\nabla{v}_j^\tau|\sigma_j)^q\}\le Lq\omega(r_j)\leq 1$ {for $j$ large enough}, then 
\[
\begin{split}
&\int_{(E_j^\tau\setminus A_h)\cap S_{j,\tau}^2}\varphi^+_j(|\nabla v_j^\tau|)\,\mathrm{d}y\\
&\le \frac{1}{\varphi(0,\sigma_j)}\int_{(E_j^\tau\setminus A_h)\cap S_{j,\tau}^2}\max\{(|\nabla{v}_j^\tau|\sigma_j)^p,(|\nabla{v}_j^\tau|\sigma_j)^q\}\varphi_{r_j}^+(1)\,\mathrm{d}y\\
&\le \gamma_d {(Lq\omega(r_j))^{\frac{p}{q}}}\frac{L}{p} \frac{1}{\varphi(0,\sigma_j)}  \underset {j\to+\infty}{\longrightarrow} 0.
\end{split}
\]
In {$S_{j,\tau}^1$ condition \vauno{} holds,}  then
\[
\begin{split}
\int_{(E_j^\tau\setminus A_h)\cap S_{j,\tau}^1}\varphi^+_j(|\nabla v_j^\tau|)\,\mathrm{d}y&\le 2\int_{(E_j^\tau\setminus A_h)\cap S_{j,\tau}^1}\varphi^-_j(|\nabla v_j^\tau|)\,\mathrm{d}y\le c\,\varphi_j^-(\tau)\,\mathcal{L}^d(E_j^\tau\setminus A_h)\\
&\le c\int_{\{M|\nabla{v}_j|>\tau\}\setminus A_h}\varphi_j^-(M|\nabla {v}_j|)\,\mathrm{d}y\,,
\end{split}
\]
where we used {\eqref{eq:liptrunc},} \eqref{ejlambda}, \eqref{eq:stimevarphij}. From the equiintegrability of the functions $\varphi_j^-(M|\nabla {\hat v}_j|)$ in ${B_R(0)\setminus A_h}$ and from \eqref{misura}, given $\eta>0$, we fix $\tau=\lambda(\eta)$ sufficiently large in order that 
\begin{equation}\label{piccolo}
c\int_{\{M|\nabla{v}_j|>\tau\}\setminus A_h}\varphi_j^-(M|\nabla {v}_j|)\,\mathrm{d}y<\eta.
\end{equation}
Therefore we can state that 
\begin{equation*}
\limsup_{j\to+\infty}\int_{E_j^\tau\setminus A_h}\varphi_j(y,|\nabla v_j^\tau|)\,\mathrm{d}y< \eta.
\end{equation*}
Concerning the first term in \eqref{liminf}, for the above fixed $\tau=\tau(\eta)$, the sequence $(v_j^\tau)_j$ is equibounded in $W^{1,\infty}(B_R(0))$, therefore, up to a subsequence, it converges to a function $v^\tau$ weakly$^*$ in $W^{1,\infty}(B_R(0))$ and in measure. Moreover, by the lower semicontinuity under convergence in measure of the map
\[
w\mapsto\mathcal{L}^d(\{x\in B_R(0)\setminus A_h : w(x)\neq 0\}),
\]
then
\begin{equation}\label{misura2}
\begin{split}
\tau^p\mathcal{L}^d(\{x\in B_1\setminus A_h : v^\tau\neq v_0\})&\le\liminf_{j\to+\infty}\tau^p\mathcal{L}^d(\{x\in B_1\setminus A_h : v_j^\tau\neq {v}_j\})\\
&\le \liminf_{j\to+\infty}\tau^p\mathcal{L}^d(E_j^\tau\setminus A_h)\\
&\le \liminf_{j\to+\infty} \frac{\tau^p}{\varphi_j^-(\tau)}\int_{\{M|\nabla {v}_j|>\tau\}\setminus A_h}\varphi_j^-(M |\nabla {v}_j|)\,\mathrm{d}y\\
&\le \liminf_{j\to+\infty}\frac{1}{\varphi_j^-(1)} \int_{\{M|\nabla {v}_j|>\tau\}\setminus A_h}\varphi_j^-(M |\nabla {v}_j|)\,\mathrm{d}y\\
&\le c \liminf_{j\to+\infty} \int_{\{M|\nabla { v}_j|>\tau\}\setminus A_h}\varphi_j^-(M |\nabla {v}_j|)\,\mathrm{d}y\le c\,\eta,
\end{split}
\end{equation}
using that $\varphi_j^-$ satisfies \inc{p}, the bound from above of $\varphi_j^-(1)$,  and \eqref{piccolo}. 
All things considered, {setting $C_s:=\{x\in B_R(0):\,\, |\nabla v_0(x)|\le s\}$,} from \eqref{liminf} we derive 
\[
\begin{split}
\liminf_{j\to+\infty}\int_{B_R(0)}\varphi_j(y,|\nabla {v}_j|&)\,\mathrm{d}y \ge \int_{B_R(0)\setminus A_h}\varphi_\infty(|\nabla v^\tau|)\,\mathrm{d}y-\eta\\
&\ge \int_{(B_R(0)\setminus A_h)\cap\{v^\tau=v_0\}\cap C_s}\varphi_\infty(|\nabla v_0|)\,\mathrm{d}y-\eta \\
&=\int_{(B_R(0)\setminus A_h)\cap C_s}\varphi_\infty(|\nabla v_0|)\,\mathrm{d}y -\int_{(B_R(0)\setminus A_h)\cap\{v^\tau\neq v_0\}\cap C_s}\varphi_\infty(|\nabla v_0|)\,\mathrm{d}y-\eta \\
& \ge \int_{(B_R(0)\setminus A_h)\cap C_s}\varphi_\infty(|\nabla v_0|)\,\mathrm{d}y - \varphi_\infty(s)\mathcal{L}^d(\{x\in B_R(0)\setminus A_h : v^\tau\neq v_0\})-\eta \\
&\ge \int_{(B_R(0)\setminus A_h)\cap C_s}\varphi_\infty(|\nabla v_0|)\,\mathrm{d}y-\varphi_\infty(s) c\,\eta-\eta,
\end{split}
\]
where we used \eqref{misura2} in the last inequality. Thus, letting first $\eta$ tend to zero, then $h$ and finally $s$ tend to infinity, we proved \eqref{lsc}. %In the same vein, \eqref{lsc} holds in every ball $B_\rho$, for $\rho\in (0,1]$, of course.
\endproof

\medskip

\paragraph{\bf{Acknowledgements}} The authors are members of Gruppo Nazionale per l'Analisi Matematica, la Probabilit\`a e le loro Applicazioni (GNAMPA) of INdAM. The work of G. Scilla and F. Solombrino is part of the project “Variational Analysis of Complex Systems in Materials Science, Physics and Biology” PRIN Project 2022HKBF5C. {The research of C. Leone and A. Verde was supported by PRIN Project 2022E9CF89 ``Geometric Evolution Problems and Shape Optimizations''.  PRIN projects are part of PNRR Italia Domani, financed by European Union through NextGenerationEU. }

%\section*{Declarations}

%\noindent {\bf  Data availability statement:} All data needed are contained in the manuscript.

%\medskip
%\noindent {\bf  Funding and/or Conflicts of interests/Competing interests:} The authors declare that there are no financial, competing or conflict of interests.

%\addcontentsline{toc}{section}{\numberline{}References}

 \bibliographystyle{siam}

\end{document}